\documentclass{amsart}
\usepackage[T1]{fontenc}
\usepackage{graphicx}
\usepackage{adjustbox, comment}

\usepackage{indentfirst}
\usepackage{mathtools}
\usepackage{tikz-cd, cancel}
\usepackage{amssymb,amsthm,amsmath,amsfonts}
\usepackage{xcolor,paralist,hyperref,fancyhdr,etoolbox}
\newtheorem{theorem}{Theorem}[section]
\newtheorem{lemma}[theorem]{Lemma}
\newtheorem{proposition}[theorem]{Proposition}
\newtheorem{corollary}[theorem]{Corollary}

\theoremstyle{definition}
\newtheorem*{definition}{Definition} 
\newtheorem{example}{Example}

\theoremstyle{remark}
\newtheorem{remark}[theorem]{Remark}

\newcommand{\PP}{\mathbb{ P}} 

\newcommand{\Hom}{\operatorname{Hom}}

\newcommand{\GL}{\operatorname{GL}}

\newcommand{\Bl}{\operatorname{Bl}}

 \newcommand{\co}{\mathcal{O}}

\makeatletter
\newcommand*{\da@rightarrow}{\mathchar"0\hexnumber@\symAMSa 4B }
\newcommand*{\da@leftarrow}{\mathchar"0\hexnumber@\symAMSa 4C }
\newcommand*{\xdashrightarrow}[2][]{%
  \mathrel{%
    \mathpalette{\da@xarrow{#1}{#2}{}\da@rightarrow{\,}{}}{}%
  }%
}
\newcommand*{\da@xarrow}[7]{%
  \sbox0{$\ifx#7\scriptstyle\scriptscriptstyle\else\scriptstyle\fi#5#1#6\m@th$}%
  \sbox2{$\ifx#7\scriptstyle\scriptscriptstyle\else\scriptstyle\fi#5#2#6\m@th$}%
  \sbox4{$#7\dabar@\m@th$}%
  \dimen@=\wd0 %
  \ifdim\wd2 >\dimen@
    \dimen@=\wd2 %
  \fi
  \count@=2 %
  \def\da@bars{\dabar@\dabar@}%
  \@whiledim\count@\wd4<\dimen@\do{%
    \advance\count@\@ne
    \expandafter\def\expandafter\da@bars\expandafter{%
      \da@bars
      \dabar@ 
    }%
  }%
  \mathrel{#3}%
  \mathrel{%
    \mathop{\da@bars}\limits
    \ifx\\#1\\%
    \else
      _{\copy0}%
    \fi
    \ifx\\#2\\%
    \else
      ^{\copy2}%
    \fi
  }%
  \mathrel{#4}%
}
\makeatother

\begin{document}
\title{Goppa Duality for Surfaces} 
\author[H. Iwasaki]{Hikari Iwasaki}
\date{\today}
\email{iwasakih@stanford.edu}
\maketitle

\let\thefootnote\relax

\begin{abstract}

 Gale duality is an involution of point configurations in projective spaces.  
Goppa duality is an extension of Gale duality, a duality between linear series on a Gorenstein curve passing through prescribed points in projective spaces. We generalize this classical result to Goppa duality for surfaces, a duality statement on linear series on surfaces which realize the prescribed points as a complete intersection of two divisors. We give a number of applications, including the existence or uniqueness of Veronese surfaces satisfying certain conditions to pass through prescribed points or curves. As a crucial example, we present an alternative proof to Coble's result on the existence of four Veronese surfaces passing nine general points in $\PP^5$.

\end{abstract} 

\bigskip

\section{Introduction}
\label{sec: 1}
\label{sec: intro}
In 1956, David Gale discovered the notion of the \emph{Gale transform} in his study of convex geometry and polytope theory \cite{Gale1956}. In algebraic terms, the Gale transform is an involution between a vector subspace and its annihilator in the dual space. Throughout this paper, we assume that $r$ and $s$ are positive integers, $\gamma = r + s + 2$, and schemes and vector spaces are defined over a fixed field $k$. If $V$ and $V'$ are $\gamma$-dimensional vector spaces with a perfect pairing $\langle \cdot,\cdot\rangle : V \times V' \to k$, so that they can be regarded as duals of each other, then for any $(r+1)$-dimensional vector subspace $W$ of $V$, its annihilator $W^\perp = \{\lambda \in V' \ | \ \langle w, \lambda \rangle = 0 \text{ for all }w \in W\}$ in $V'$ is $(s+1)$-dimensional, and they fit in the short exact sequence
$$0 \to W \to V \cong (V')^\vee \to (W^\perp)^\vee \to 0,$$
where the second map is the dual of the inclusion of $W^\perp$ in $V'$. The linear maps $W \to V$ and $W^\perp \to V'$ define configurations of $\gamma$ points in $\PP W^\vee = \PP^r$ and $\PP(W^\perp)^\vee = \PP^s$, respectively. This involution on $\gamma$ points in the projective spaces is called \emph{Gale duality}.

In 1981, Valery Denisovich Goppa discovered an algebro-geometric application of Gale duality in terms of divisors on curves, which is now called \emph{Goppa duality} \cite{Goppa1981}.   Roughly speaking, Goppa duality can be interpreted as an extension of Gale duality, from a duality of maps from points to projective spaces to a duality of maps from curves to projective spaces. 
Specifically, if we have a closed embedding of a zero-dimensional Gorenstein scheme $\Gamma$ into a Gorenstein curve $C$ and $L$ is a line bundle on $C$, then the long exact sequence in cohomology associated to the short exact sequence for the closed embedding of $\Gamma$ in $C$, together with Serre duality, induces the following exact sequence:
\begin{align*}
    H^0(C, L) \xrightarrow{\operatorname{res}} H^0(\Gamma, L|_\Gamma) 
    \cong H^0(\Gamma, K_\Gamma \otimes L|_\Gamma^\vee)^\vee
    \xrightarrow{\operatorname{res}^\vee} H^0(C, K_C(\Gamma)\otimes L^\vee)^\vee .
\end{align*}
Restating the exactness of the sequence, we obtain a pair of a linear series $W$ of $L$ and a linear series $W^\perp$ of $K_C(\Gamma) \otimes L$ complementary to the kernel of the corresponding restriction map of global sections. It follows that $W$ and $W^\perp$ fit into a short exact sequence as in the second row of the following diagram:

\adjustbox{scale=0.9, center}{
\begin{tikzcd}[column sep=small]
{H^0(C, L(-\Gamma))} \arrow[r] 
  & {H^0(C,L)} \arrow[r, "\operatorname{res}"] 
  & {H^0(\Gamma, L|_\Gamma) \, \cong \, H^0(\Gamma, K_\Gamma \otimes L|_\Gamma^\vee)^\vee} 
      \arrow[r, "\operatorname{res}^\vee"] \arrow[d, equal] 
  & {H^0(C, K_C(\Gamma) \otimes L^\vee)^\vee} \arrow[r] \arrow[d, two heads] 
  & {H^0(C, K_C \otimes L^\vee)^\vee} \\
0 \arrow[r] 
  & W \arrow[u, hook] \arrow[r, "\operatorname{res}"] 
  & {H^0(\Gamma, L|_\Gamma) \, \cong \, H^0(\Gamma, K_\Gamma \otimes L|_\Gamma^\vee)^\vee} 
      \arrow[r, "\operatorname{res}^\vee"] 
  & {(W^\perp)^\vee} \arrow[r] 
  & 0
\end{tikzcd}
}
\bigskip

\noindent In particular, the short exact sequence induces a Gale dual pair of maps $\Gamma \to \PP W^\vee$ and $\Gamma \to \PP(W^\perp)^\vee$,  which factors through $C$:
$\Gamma  \xhookrightarrow{\alpha} C \xdashrightarrow{W} \PP W^\vee$ and $\Gamma  \xhookrightarrow{\alpha} C \xdashrightarrow{W^\perp}\PP (W^\perp)^\vee$. This duality between the two maps from $C$ to the projective spaces is called (classical) \emph{Goppa duality}. 
This discussion will be expanded in detail in Section \ref{sec: 2}.

As introduced above, classical Goppa duality is concerned with factorizations of maps from $\Gamma$ to projective spaces through a curve $C$.
In Section \ref{sec: 3}, we present the main result of this paper, \emph{Goppa duality for surfaces}, which extends classical Goppa duality to factorizations of maps through a surface $S$. More precisely, we study factorizations of maps $\Gamma \to \PP^r$ through a Gorenstein surface $S$, specifically when $\Gamma$ is a complete intersection of two curves on $S$. This leads to a duality between two linear series $W$ and $W^\perp$ on $S$, and hence a duality between the associated maps $\Gamma  \xhookrightarrow{\alpha} S \xdashrightarrow{W} \PP W^\vee$ and $\Gamma  \xhookrightarrow{\alpha} S \xdashrightarrow{W^\perp}\PP (W^\perp)^\vee$.

To illustrate usefulness of this extension, we focus on the case where $S$ is the projective plane $\mathbb{P}^2$ and derive several geometric results. These include:
\begin{itemize}
    \item 
    (Corollary \ref{cor: unique Veronese through two deg 6 elliptic curves in P5}) If two nondegenerate sextic elliptic curves in $\PP^5$ intersect at exactly nine nondegenerate points in $\PP^5$, then there exists a unique Veronese surface containing the elliptic curves.   
    
    \item 
    (Proposition \ref{prop: four veronese through 9 points}) An alternative proof to \cite{Coble1922}: exactly four Veronese surfaces pass through nine general points in $\mathbb{P}^5$. 
    
    \item 
    (Proposition \ref{prop: existence of factorization through BlpP2}) Eight general points in $\mathbb{P}^4$ admit a unique factorization through a surface isomorphic to $\mathrm{Bl}_p \mathbb{P}^2$ via the complete linear series $|2H-E|$, where the points lie on a complete intersection of two curves of class $3H - E$.
    
\end{itemize}

The paper is structured as follows. In Section  \ref{sec: 2}, we review Gale duality and classical Goppa duality with explicit examples; a deep analysis on these will motivate and set up later arguments for extending the duality. Section  \ref{sec: 3} describes the main result of the paper: a generalization of Goppa duality to surfaces. To derive the result, in Section \ref{sec: 3.1}, we first describe an algebraic result on cohomology proved by David Eisenbud and Sorin Popescu \cite{EisenbudPopescu2000}, which can be used to extend Goppa duality to higher dimensions. Section \ref{sec: complete intersection} specializes the algebraic result in Section \ref{sec: 3.1} to the case where the finite Gorenstein scheme $\Gamma$ is realized as a complete intersection on a Gorenstein surface $S$, leading to the main result of this paper, {Goppa duality for surfaces}. We give detailed applications to the case when $S$ is the projective plane $\PP^2$ in Section \ref{sec: S=PP2}.  Finally, Section \ref{sec: almost complete intersection} addresses the more complex case where the points $\Gamma$ in $\PP^2$ are ``almost'' a complete intersection, in the sense that they become a complete intersection when we add excess point $R$. By reframing the implication of Goppa duality for surfaces, we deduce an alternative proof to Arthur B. Coble's result  \cite{Coble1922} on the existence of four Veronese surfaces through nine general points in $\PP^5$.
Finally, we conclude the paper with Section \ref{sec: 4}, where we present an interpretation of a prior work by Anand Deopurkar and Anand Patel \cite{DeopurkarPatel2024} on the count of Veronese surfaces passing through thirteen general points in $\PP^5$ using Goppa duality for surfaces, as well as potential further generalizations of the theory. 

\subsection*{Acknowledgements}
I am grateful to my advisor Ravi Vakil for suggesting a related topic and for helpful discussions.

\section{Gale Transform and Goppa Duality For Curves}
\label{sec: 2}
\subsection{Gale Duality} 
\label{sec: Gale Transform}

Gale duality is a duality of point configurations in two projective spaces, which can be motivated and explicitly computed in terms of linear algebra as follows. 
Recall from Section \ref{sec: intro} that $r$ and $s$ are positive integers and $\gamma = r + s + 2$. Suppose that we have $\gamma$ non-degenerate points $\{p_i\}_{i=1}^\gamma$ in $\PP^r_k$, with coordinates $p_i = [x_{i0}: \cdots : x_{ir}]$. Then we can construct a $\gamma \times (r+1)$ matrix $G$ whose rows are $\begin{bmatrix}
    x_{i0} &  \cdots & x_{ir}
\end{bmatrix} $. By nondegeneracy, $G$ has trivial kernel, hence by the Rank-Nullity Theorem, the cokernel map can be represented by an $(s+1) \times \gamma$-matrix; denote the transpose of this matrix by $G'$. It follows that we have a short exact sequence
$$0 \to k^{\oplus (r+1)} \xrightarrow{G} k^{\oplus \gamma} \xrightarrow{(G')^\top} k^{\oplus (s + 1)} \to 0.$$ 
The matrix $G'$ induces $\gamma$ points $\{q_i\}_{i=1}^\gamma$ on $\PP^{s}$ whose coordinates $[y_{i0}: \cdots : y_{is}]$ are the rows $\begin{bmatrix}
    y_{i0} &  \cdots & y_{is}
\end{bmatrix} $ of $G'$. 

More generally, observe that scaling the columns of $G$ or $G'$ does not affect the exactness of the sequence, because such an operation leaves the image of the map invariant. Hence, we have $(G')^\top D G = 0$ for all invertible diagonal matrix $D$, and conversely, any short exact sequence 
$0 \to k^{\oplus (r+1)} \xrightarrow{} k^{\oplus \gamma} \xrightarrow{} k^{\oplus (s+1)} \to 0$ would induce a pair $(G,G')$ of two matrices, whose rows induce $\gamma$ points in $\PP^r$ and $\PP^s$, respectively.
We will call this pair of point configurations $\Gamma = \{[x_{i0}:\cdots : x_{ir}]\}_{i=1}^\gamma$ and $\Gamma' = \{[x_{i0}:\cdots : x_{is}]\}_{i=1}^\gamma$ the \textit{Gale transform} of each other.

\begin{definition}
    Let $\Gamma \subset \PP^r$ and $\Gamma' \subset \PP^s$ be two nondegenerate configurations of $\gamma$  points in the projective spaces. Suppose that $G $ and $ G'$ are $\gamma \times (r+1)$ and $\gamma \times (s+1)$-matrices, whose rows describe homogeneous coordinates of points  $\Gamma$ and $\Gamma'$, respectively. Then $\Gamma'$ is a \emph{Gale transform} of $\Gamma$ if $(G')^\top DG = 0$ for some invertible diagonal matrix $D$. 
\end{definition}

It follows that any nondegenerate $\gamma$ points in $\PP^r$  have a Gale transform in $\PP^s$, and this is a well-defined duality notion once we impose projective equivalence on the projective spaces, as proved in the following proposition. 

\begin{proposition}
    The Gale transform of $\Gamma \subset \PP^r$ is well-defined and unique up to projective linear transformation on $\PP^r$ and $\PP^s$.
\end{proposition}

\begin{proof}
    If rows of $G_1$ and $G_2$ describe the homogeneous coordinates of the same point configuration $\Gamma$ in $\PP^r$, then $G_1 = DG_2$ for some invertible diagonal matrix $D$. Thus $(G_1')^\top D' G = 0$ for some invertible diagonal matrix $D'$ if and only if $(G_2')^\top D' G = 0$ for some invertible diagonal matrix $D'$, and the same argument holds for $G'$. Therefore, the Gale transform is well-defined. Moreover transformation of $G$ or $G'$ by invertible matrices in $\GL_{r+1}$ or $\GL_{s+1}$ acting from the right leaves the condition invariant up to projective linear transformations. 
\end{proof}

\begin{example}\label{ex: gale transform of 5 points}
    Consider the set of five points 
    $$\Gamma= \{[1:0:0], [0:1:0], [0:0:1], [1:1:1], [1:a:b]\}$$ in $\PP^2$. This induces a $5 \times 3$-matrix
    $$\begin{bmatrix}
        1 & 0 & 0 \\
        0 & 1 & 0 \\
        0 & 0 & 1 \\
        1 & 1 & 1 \\
        1 & a & b 
    \end{bmatrix}, $$\text{and the kernel of its transpose is the column space of }$$ \begin{bmatrix}
        -1 & -1 \\
        -1 & -a \\
        -1 & -b \\
        1 & 0 \\
        0 & 1  
    \end{bmatrix}.$$
    Therefore, the Gale transform of $\Gamma$ is the set $\Gamma'$ of five points in $\PP^1$, where
    $$\Gamma' = \{[1:1], [1:a], [1:b], [1:0], [0:1]\}.$$
\end{example}

\bigskip

The uniqueness of the Gale transform up to projective linear transformation suggests that Gale duality is determined by the images of $G$ and $G'$ in the $\gamma$-dimensional vector space $k^{\oplus \gamma}$.  Therefore, we may understand Gale duality in a coordinate-free form, which leads to the following algebraic interpretation of the duality as the pair of a vector subspace and its annihilator under a perfect pairing. The proof is an algebraic exercise left to the reader.

\begin{lemma}
    \label{lem: coordinate-free Gale duality}
    Suppose that $V$ is a $\gamma$-dimensional vector space, and $V'$ its dual space under a perfect pairing $\langle \cdot, \cdot \rangle$. If $W$ is an $(r+1)$-dimensional vector subspace of $V$ and the subspace $W^\perp$ of $V'$ is the $(s+1)$-dimensional annihilator of $W$, i.e., $W^\perp = \{\lambda \in V' \ | \ \langle \lambda, w\rangle = 0 \text{ for all }w \in W\}$, then the pair $(W, W')$ satisfies the short exact sequence
    $$0 \to W \hookrightarrow V \cong (V')^\vee \to (W^\perp)^\vee \to 0,$$
    where the second map is the dual of the inclusion map of $W^\perp$ into $V'$.
\end{lemma}

This yields a useful interpretation of Gale duality: it is the data of a pair $(W, W^\perp)$ of subspaces of $V$ and $V'$ with a perfect pairing which are annihilators of each other. The point configurations in projective spaces can then be understood as the images of the standard coordinates $\Gamma = \{[1:0:\cdots : 0], \cdots, [0:\cdots : 0 : 1]\}$ under the induced rational maps $\PP V^\vee \dashrightarrow \PP W^\vee$  and $\PP (V')^\vee \dashrightarrow \PP (W^\perp)^\vee$  of projective spaces.

\subsection{Gale duality in algebro-geometric terms}

To extend Gale duality to Goppa duality, it helps to understand Gale duality in algebro-geometric terms first.
From the discussion in the previous section, Gale duality can be rephrased as the pair of two maps from the set $\Gamma$ of $\gamma$ $k$-rational points into projective spaces: 
\[\begin{tikzcd}[column sep=huge, row sep=small]
\PP^r = \PP W^\vee &  &  \PP (W^\perp)^\vee  = \PP^s\\
      &  &             \\
      & \Gamma  \arrow[luu, " ", rounded corners=10pt, 
    to path={(\tikztostart) -| (\tikztotarget) \tikztonodes}] \arrow[ruu, " "', rounded corners=10pt, 
    to path={(\tikztostart) -| (\tikztotarget) \tikztonodes}] &     
\end{tikzcd}
\]
On the other hand, recall that maps to projective spaces can be determined by the data of a line bundle and sections. Under this perspective, the two maps $\Gamma \to \PP W^\vee$ and $\Gamma \to \PP (W^\perp)^\vee$ correspond to the data $(L, \sigma_0,\cdots, \sigma_r)$ and  $(L', \tau_0,\cdots, \tau_s)$, respectively,  where $L$ and $L'$ are line bundles on $\Gamma$, and  $\{\sigma_i\}_{i=1}^r$ and  $\{\tau_i\}_{i=1}^s$ respectively form bases for the linear series $W$ and $W^\perp$. To fit this back into the perspective of short exact sequences given in Lemma \ref{lem: coordinate-free Gale duality}, we take $V$ and $V'$ to be the spaces of global sections $H^0(\Gamma, L)$ and $H^0(\Gamma, L')$, respectively. Gale duality then requires a perfect pairing between $H^0(\Gamma, L)$ and $H^0(\Gamma, L')$; we will use the most famous duality in algebraic geometry: Serre duality. Because $\Gamma$ consists of $k$-rational (smooth) points, we have a perfect pairing $H^0(\Gamma, L) \times  H^0(\Gamma, K_\Gamma \otimes L^\vee) \to k$. Because all line bundles on the finite scheme $\Gamma$ of $k$-rational points are isomorphic to the trivial line bundle, we may set $L'$ to be $K_\Gamma \otimes L^\vee$.  Under this identification, we obtain a short exact sequence
$$ 0 \to W \to H^0(\Gamma, L) \cong H^0(\Gamma, K_\Gamma \otimes L^\vee)^\vee \to (W^\perp)^\vee \to 0. $$
Hence in algebraic geometry, we may understand Gale duality as the pair $(W, W^\perp)$ of linear series of $L$ and $K_\Gamma \otimes L^\vee$ on $\Gamma$ under the perfect pairing induced by the Serre duality on $\Gamma$. We can decorate the diagram with the information of linear series inducing the maps:
\[\begin{tikzcd}[column sep=huge, row sep=small]
\PP W^\vee &  & \PP (W^\perp)^\vee \\
      &  &             \\
      & \Gamma  \arrow[luu, "W \subset |L|", rounded corners=10pt, 
    to path={(\tikztostart) -| (\tikztotarget) \tikztonodes}] \arrow[ruu, "W^\perp \subset |K_\Gamma \otimes L^\vee|"', rounded corners=10pt, 
    to path={(\tikztostart) -| (\tikztotarget) \tikztonodes}] & 
\end{tikzcd}
\]

\bigskip
Finally, because the key component of this argument is the validity of Serre duality on $\Gamma$, the algebro-geometric Gale duality applies not only to the set $\Gamma$ of $\gamma$ closed points with residue field $k$, but also to any finite Gorenstein $k$-scheme, on which as Serre duality holds on all line bundles. In summary, we obtain the algebro-geometric version of Gale duality. 

\begin{theorem}[Gale duality, algebro-geometric version]
    Gale duality is a duality between two maps from a Gorenstein finite scheme $\Gamma$ to projective spaces,  $$\alpha: \Gamma \xrightarrow{W} \PP^r \text{ and } \beta : \Gamma \xrightarrow{W^\perp} \PP^s,$$
    where $W$ is a linear series of $L$, $W^\perp$ is the linear series of $K_\Gamma \otimes L^\vee$, and $W$ and $W^\perp$ are annihilators of each other under Serre duality.
\end{theorem}

\subsection{Classical Goppa Duality (Goppa Duality for Curves)}
\label{sec: Goppa duality}

In the previous section, we introduced Gale duality as a duality between maps of Gorenstein points to projective spaces. In this section, we derive that Gale duality can be promoted to a duality between maps  from a Gorenstein curve to projective spaces, called (classical) \emph{Goppa duality}.


Suppose that there exists a closed embedding of a finite (zero-dimensional) Gorenstein scheme $\Gamma$ into a Gorenstein curve $C$, where $\Gamma$ is an effective Cartier divisor of $C$. Let $L$ be a line bundle on $C$. The closed subscheme exact sequence induces the long exact sequence of cohomology,
\begin{align*}
    0 &\to H^0(C, L(-\Gamma)) \to H^0(C, L)  \xrightarrow{\operatorname{res}} H^0(\Gamma, L|_\Gamma)  
    \to H^1(C,L(-\Gamma)) \to H^1(C, L)  \to 0.
\end{align*}
By Serre duality, this reduces to

\medskip

\adjustbox{scale=0.9, center}{%
\begin{tikzcd}[column sep=small]
{0} 
  \arrow[r] & {H^0(C,L(-\Gamma))} 
  \arrow[r] & {H^0(C, L)} 
  \arrow[r, "\operatorname{res}"] & {H^0(\Gamma, L|_\Gamma) \cong H^0(\Gamma, K_\Gamma \otimes L|_\Gamma^\vee)^\vee} 
  \arrow[r, "\operatorname{res}^\vee"] & {H^0(C, K_C(\Gamma)\otimes L^\vee)^\vee} 
  \arrow[r] & {H^0(C, K_C \otimes L^\vee)^\vee} 
  \arrow[r] & {0}
\end{tikzcd}%
}
\smallskip

\noindent where the maps on the dual spaces are given by the dual of the maps
$$H^0(C, K_C \otimes L^\vee) \to H^0(C, K_C(\Gamma)\otimes L^\vee) \xrightarrow{\operatorname{res}} H^0(\Gamma, K_\Gamma \otimes L|_\Gamma^\vee).$$

Using the middle part of the exact sequence, we can understand Gale transform of points in terms of curves passing through them. More precisely, let $W \subset H^0(C,L)$ be a linear series  complementary to the kernel of the restriction map $H^0(C,L) \to H^0(\Gamma, L|_\Gamma)$,  and similarly let $W^\perp \subset H^0(C,K_C(\Gamma)\otimes L^\vee)$ be a linear series complementary to the kernel of the restriction map $H^0(C,K_C(\Gamma) \otimes L^\vee) \to H^0(\Gamma, K_\Gamma \otimes L|_\Gamma^\vee)$. Then $W$ and $W^\perp$ fit in a short exact sequence which induce Gale duality, given by the second row of the following diagram:

\bigskip

\adjustbox{scale=0.9, center}{
\begin{tikzcd}[column sep=small]
{H^0(C, L(-\Gamma))} \arrow[r] & {H^0(C,L)} \arrow[r, "\operatorname{res}"]        & {H^0(\Gamma, L|_\Gamma) \cong H^0(\Gamma, K_\Gamma \otimes L|_\Gamma^\vee)^\vee} \arrow[r, "\operatorname{res}^\vee"] \arrow[d, equal] & {H^0(C,K_C(\Gamma) \otimes L^\vee)^\vee} \arrow[r] \arrow[d, two heads] & {H^0(C,K_C \otimes L^\vee)^\vee} \\
0 \arrow[r]                    & W \arrow[u, hook] \arrow[r, "\operatorname{res}"] & {H^0(\Gamma, L|_\Gamma) \cong H^0(\Gamma, K_\Gamma \otimes L|_\Gamma^\vee)^\vee} \arrow[r, "\operatorname{res}^\vee"]                                & (W^\perp)^\vee \arrow[r]                                                       & 0                            .   
\end{tikzcd}
}
\bigskip

Then the induced Gale dual pair $\Gamma \to \PP W^\vee$ and $\Gamma \to \PP (W^\perp)^\vee$ has the following geometric interpretation: the two maps are the restriction to $\Gamma$ of the rational maps $C \xdashrightarrow{W} \PP W^\vee$ induced by the linear series $W$ and  $C \xdashrightarrow{W^\perp}  \PP (W^\perp)^\vee$ induced by the linear series $W^\perp$. In other words, the Gale dual pair factors through the curve $C$ via linear series $W$ and $W^\perp$.  
This result is called \emph{Goppa duality,} which we refer to also as \emph{classical Goppa duality} or \emph{Goppa duality for curves} in this paper. 

\begin{theorem} [classical Goppa duality]
    \label{thm: Goppa duality for curves}
    Let $\iota: \Gamma \to C$ be an effective Cartier divisor of length $\gamma$ into a Gorenstein curve $C$. Suppose that $W$ is an $(r+1)$-dimensional linear series of $L$ complementary to the kernel of the restriction map $H^0(C,L) \to H^0(\Gamma, L|_\Gamma)$. Then there exists a $(s+1)$-dimensional linear series $W^\perp$ of $K_C (\Gamma) \otimes L^\vee$ satisfying the following short exact sequence,
    \begin{align*}
        \label{tag: ast}
        0 \to W \xrightarrow{\operatorname{res}} H^0(\Gamma, L|_\Gamma) \cong H^0(\Gamma, K_\Gamma\otimes L|_\Gamma^\vee)^\vee \xrightarrow{\operatorname{res}^\vee} (W^\perp)^\vee \to 0, \tag{$\ast$}
    \end{align*}
    and such a linear series $W^\perp$ is complementary to the kernel of  the restriction map $H^0(C,K_C(\Gamma)\otimes L^\vee) \to H^0(\Gamma, K_\Gamma \otimes L|_\Gamma^\vee)$.
    
    In particular, the Gale dual pair $\Gamma \to \PP^r$ and $\Gamma \to \PP^s$ induced by \emph{($\ast$)} factors through $C$ as 
    $$\Gamma \xrightarrow{\iota} C \xdashrightarrow{W} \PP^r \text{ and } \Gamma \xrightarrow{\iota} C \xdashrightarrow{W^\perp} \PP^s, $$
    where the rational maps $C \xdashrightarrow{W} \PP^r$ and $ C \xdashrightarrow{W^\perp} \PP^s$ are induced by the linear series $W$ and $W^\perp$, respectively.

\end{theorem}

Note that the theorem statement implicitly uses the fact that an effective Cartier divisor of a Gorenstein scheme, or more generally a regular embedding of a Gorenstein scheme, is automatically Gorenstein \cite[Section 21.3]{Eisenbud}, to guarantee that $\Gamma$ is a Gorenstein scheme. 

\medskip

\begin{proof}
    The closed subscheme exact sequence for $\iota : \Gamma \to C$  along with Serre duality induces the following exact sequence of cohomology,
    \[ H^0(C, L)  \xrightarrow{\operatorname{res}} H^0(\Gamma, L|_\Gamma) \cong H^0(\Gamma, K_\Gamma \otimes L|_\Gamma^\vee)^\vee 
    \xrightarrow{\operatorname{res^\vee}} H^0(C, K_C(\Gamma)\otimes L^\vee)^\vee .
    \]

    By the assumption on $W$, $W$ is isomorphic to the image of the restriction map $H^0(C,L) \to H^0(\Gamma, L|_\Gamma)$.
On the other hand, choose $W^\perp$ to be any linear series $K_C(\Gamma) \otimes L^\vee$ complementary to the kernel of the restriction map $H^0(C,K_C(\Gamma) \otimes L^\vee) \to H^0(\Gamma, K_\Gamma \otimes L|_\Gamma^\vee)$. Then $W^\perp$ is isomorphic to the image of the restriction map, so the dual of the restriction map factors through $(W^\perp)^\vee$.
    It follows that $W$ and $W^\perp$ fit into the short exact sequence induced by the exactness of (\ref{tag: ast}) at $H^0(\Gamma, L|_\Gamma) \cong H^0(\Gamma, K_\Gamma \otimes L|_\Gamma^\vee)^\vee$ . Conversely, any linear series $W^\perp$ of $K_G(\Gamma) \otimes L^\vee$ fitting into (\ref{tag: ast}) must be isomorphic to the image of the restriction map $H^0(C,K_C(\Gamma) \otimes L^\vee) \to H^0(\Gamma, K_\Gamma \otimes L|_\Gamma^\vee)$, which is equivalent to being complementary to the kernel of the restriction map.

    Next, from the short exact sequence (\ref{tag: ast}),
    we have an induced Gale dual pair $\Gamma \to \PP^r$ and $\Gamma \to \PP^s$. On the other hand, the maps $C \xdashrightarrow{W} \PP^r$ and $C \xdashrightarrow{W^\perp} \PP^s$ induced by the linear series $W$ and $W^\perp$, when restricted to $\Gamma$ via $\iota$, recover this Gale dual pair; therefore, the Gale dual pair factors through $C$ as claimed.
\end{proof}
In other words, Goppa duality for curves is a duality between two maps $$C \xdashrightarrow{W} \PP^r \text{ and }  C \xdashrightarrow{W^\perp} \PP^s,$$
where $W$ is a linear series of $L$ and $W^\perp$ is the linear series of $K_C(\Gamma) \otimes L^\vee$, whose restrictions to $\Gamma$ induce an annihilator pair under Serre duality on $\Gamma$. We may understand the relation of  Goppa duality to Gale duality as filling of the level diagram presented below.

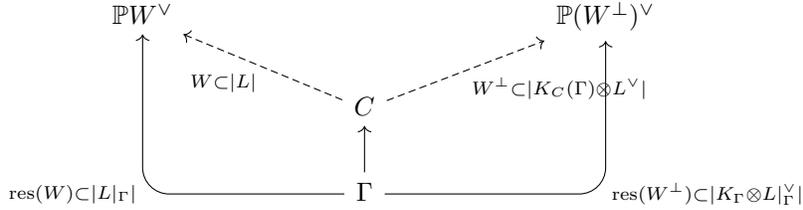
\[\begin{tikzcd}[column sep=6em]
\PP W^\vee &  & \PP (W^\perp)^\vee \\
      & C \arrow[lu, dashed, "W \subset |L|"] \arrow[ru, dashed, "W^\perp \subset |K_C(\Gamma) \otimes L^\vee|"'] &             \\
      & \Gamma \arrow[u] \arrow[luu, "\operatorname{res}(W) \subset |L|_\Gamma|", rounded corners=10pt, 
    to path={(\tikztostart) -| (\tikztotarget) \tikztonodes}] \arrow[ruu, "\operatorname{res}(W^\perp) \subset |K_\Gamma \otimes L|_\Gamma^\vee|"', rounded corners=10pt, 
    to path={(\tikztostart) -| (\tikztotarget) \tikztonodes}] & 
\end{tikzcd}
\]

\bigskip

\begin{remark}
    \label{rmk: kernel is twist by -Gamma}
    From the long exact sequence of cohomology used in the proof, the kernel of the restriction map $ H^0(C, L)  \xrightarrow{\operatorname{res}} H^0(\Gamma, L|_\Gamma)$ is exactly $H^0(C, L(\Gamma))$, and similarly the kernel of the restriction map $ H^0(C, K_C(\Gamma)\otimes L^\vee)  \xrightarrow{\operatorname{res}} H^0(\Gamma, K_\Gamma \otimes L|_\Gamma^\vee)$ is exactly $H^0(C, K_C\otimes L^\vee)$. In particular, the two restriction maps are injective if $L$ and $K_\Gamma \otimes L|_\Gamma^\vee$ have negative degrees, and in such a case, $W$ and $W^\perp$ must be the complete linear series $|L|$ and $|K_C(\Gamma) \otimes L^\vee|$, respectively, for Goppa duality (Theorem \ref{thm: Goppa duality for curves}) to be applied.
\end{remark}

\begin{remark}
\label{rmk: curves}
We emphasize that the choice of $W^\perp$ is not unique in general, especially when the restriction map of the global sections have nontrivial kernel; in such a case, we need to make a choice for a complementary subspace of the kernel. On the other hand, we note that the linear series is indeed unique when the corresponding restriction map is injective. In particular, by Remark \ref{rmk: kernel is twist by -Gamma}, when $L$ and $K_C(\Gamma)\otimes L^\vee$ have negative degrees, Goppa duality induces a one-to-one correspondence between the two maps $C \xdashrightarrow{W} \PP^r$ and  $C \xdashrightarrow{W^\perp} \PP^s$ factoring the Gale dual pair. This observation will be particularly important in Proposition \ref{prop: uniqueness of RNC} when we derive uniqueness of the factorization. 
\end{remark}

In practice, we use Goppa duality for curves (Theorem \ref{thm: Goppa duality for curves}) as follows: given $\gamma$ points $\Gamma \to \PP^r$ in a projective space, we try to find a Gorenstein curve $C$ with the following properties:
\begin{enumerate}[(i)]
    \item $C$ factors the map $\Gamma \to \PP^r$ via a closed embedding $\iota: \Gamma \hookrightarrow C$, through which $\Gamma$ can be regarded as an effective Cartier divisor, and
    \item the map $C \dashrightarrow \PP^r$ in the factorization is induced by an $(r+1)$-dimensional linear series $W$ of some line bundle $L$, complementary to the subspace $H^0(C,L(-\Gamma))$ in $H^0(C,L)$. 
\end{enumerate} 
Once we find such a Gorenstein curve, then we can factor the Gale transform $\Gamma \to \PP^s$ via $C$ as the composition of the closed embedding $\iota$ and the map $C \dashrightarrow \PP^s$, which is induced by an $(s+1)$-dimensional linear series $W^\perp$ of $K_C(\Gamma) \otimes L^\vee$ complementary to the subspace $H^0(C, K_C\otimes L^\vee) $ in $ H^0(C, K_C(\Gamma) \otimes L^\vee)$.
In particular, if the problem of finding a Gorenstein curve satisfying the two conditions (i) and (ii) is relatively easy to solve on one side of the duality, then we can answer questions on the factorization of points through curves on the other side of the duality.  

\begin{example} 
    \label{ex: goppa duality of 5 points}
    Let's return to Example \ref{ex: gale transform of 5 points} on the five points 
    $$\Gamma= \{[1:0:0], [0:1:0], [0:0:1], [1:1:1], [1:a:b]\} \text{ in }\PP^2.$$  We observed that the Gale transform of $\Gamma$ is 
    $$\Gamma' = \{[1:1], [1:a], [1:b], [1:0], [0:1]\} \text{ in } \PP^1.$$ 

    We analyze if and how this Gale dual pair factors through some Gorenstein curve $C$. We note that on one side of the duality, we have $\Gamma$ mapping into the projective line $\PP^1$. In this setting, we can regard the projective line as the Gorenstein curve, with the map $C \to \PP^1$ being the identity map and the closed embedding $\Gamma \hookrightarrow C$ being the given map itself. Then condition (i) holds because any finite set of points on $\PP^1$ form an effective Cartier divisor. Next, we note that the identity map $C \to \PP^1$ is a map induced by the complete linear series $|\co(1)|$. Because $\co(1)(-\Gamma)$ has a negative degree and $W$ is the complete linear series $|\co(1)|$, condition (ii) is also satisfied. 
    Therefore, by Goppa duality for curves (Theorem \ref{thm: Goppa duality for curves}) with the choice of the curve being $\PP^1$ and a factorization of the map to $\PP^1$ through the curve by $\Gamma \xrightarrow{\iota} \PP^1 \xrightarrow{|\co(1)|} \PP^1$, the Gale transform $\Gamma \to \PP^2$ factors through the curve $C=\PP^1$. The map $C \to \PP^2$ in the factorization of the Gale transform is induced by a 3-dimensional complete linear series of $K_{\PP^1}(\Gamma) \otimes \co(1)^\vee = \co_{\PP^1}(2)$, which must be the complete linear series $|\co_{\PP^1}(2)|$ by the dimension requirement. Therefore, the Gale transform factors as $\Gamma \xrightarrow{\iota} \PP^1 \xrightarrow{|\co(2)|} \PP^2$, from which we see that the second map
    $\PP^1 \to \PP^2$ is the rational normal curve. In other words, Goppa duality for curves deduces that there exists a conic through five nondegenerate points in $\PP^2$, and the position of the points on the conic regarded as $\PP^1$ is determined by the Gale transform up to projective equivalence. The following level diagram summarizes the key argument of the duality.

    
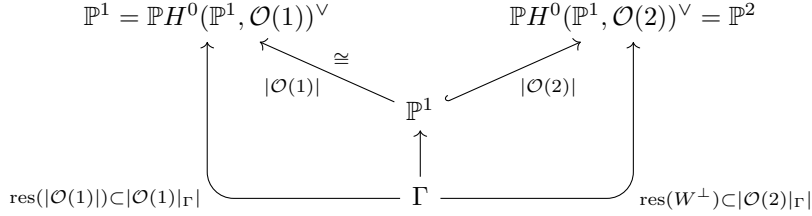
\[\begin{tikzcd}[column sep=2em]
\PP^1 = \PP H^0(\PP^1, \co(1))^\vee &  & \PP H^0(\PP^1, \co(2))^\vee = \PP^2 \\
      & \PP^1 \arrow[lu, "|\co(1)|", "\cong"'] \arrow[ru, hook, "|\co(2)|"'] &             \\
      & \Gamma \arrow[u] \arrow[luu, "\operatorname{res}(|\co(1)|) \subset |\co(1)|_\Gamma|", rounded corners=10pt, 
    to path={(\tikztostart) -| (\tikztotarget) \tikztonodes}] \arrow[ruu, "\operatorname{res}(W^\perp) \subset |\co(2)|_\Gamma|"', rounded corners=10pt, 
    to path={(\tikztostart) -| (\tikztotarget) \tikztonodes}] & 
\end{tikzcd}
\]

    \vspace{1em}
    
    We can confirm this statement also by explicit computation. Recall that any five general points in $\PP^2$ determine a unique conic. The unique conic passing through $\Gamma$ is given by 
    \begin{align*}
        b(a-1)X(Y-Z) + (a-b) Z (X-Y) = 0.
    \end{align*}
    Restricting to the case where the conic is smooth and regarding the conic as a rational normal curve in $\PP^2$, explicit computation shows that the parametrization of $\PP^1 \to \PP^2$ sending the standard coordinates $[1:0],[0:1],[1:1]$ to $[1:1:1],[1:a:b],[1:0:0]$, respectively, and passing through $[0:1:0]$ and $[0:0:1]$, is given uniquely by 
    \begin{align*}
        \PP^1 &\to \PP^2 \\
        [s:t] &\mapsto 
        [(as - t)(bs - t) : a(bs -t) (s-t) : b (as - t)(s-t)].
    \end{align*}
    This map sends each point of $\Gamma' \subset \PP^1$ to $\Gamma \subset \PP^2$, as expected by Goppa duality: 
    \begin{align*}
        [1:1] &\mapsto [1:0:0], \\
        [1:a] &\mapsto [0:1:0], \\
        [1:b] &\mapsto [0:0:1], \\
        [1:0] &\mapsto [1:1:1], \\
        [0:1] &\mapsto [1:a:b].
    \end{align*}
\end{example}

\medskip 

Although an explicit construction of the Gale transform is possible as above, it is algebraically cumbersome, and the computation does not induce a clear geometric interpretation of the result. This highlights the strength of Goppa duality in contrast: Goppa duality for curves allows us to understand the duality of points in algebro-geometric terms, involving linear series on curves and the induced maps to projective spaces.

\bigskip

Note that the key idea presented in Example \ref{ex: goppa duality of 5 points} for the construction of a rational normal curve is to identify the target projective space of the Gale transform as the curve itself, which restricts the curve to be $\PP^1$ and $r$ to be 1. In particular, we can generalize the result of Example  \ref{ex: goppa duality of 5 points} to projective spaces of higher dimension, leading to the following proof of the classical result regarding existence of a rational normal curve in $\PP^s$ passing through $s + 3$ general points. 
\begin{lemma}  
    \label{lem: existence of RNC}
    Given $s + 3$ general points in $\PP^s$, there exists a rational normal curve passing through them.
\end{lemma}

\begin{proof}
    Let $\beta : \Gamma \hookrightarrow \PP^s$ be the closed embedding of $s+3$ general points $\Gamma$ in $\PP^s$ and denote its Gale transform by $\alpha: \Gamma \to \PP^1$. By taking the curve to be $\PP^1$ and the morphism to $\PP^1$ to be the identity map, we can set 
    $$L = \co(1), \ W = H^0(\PP^1, \co(1)), \ K_C(\Gamma)\otimes L^\vee = \co(-2)\otimes\co(s+3)\otimes\co(-1)=\co(s),$$
    where $W$ satisfies the requirement stated in Goppa duality.
    By Goppa duality for curves, $W^\perp$ is an $(s+1)$-dimensional subspace of $H^0(\PP^1, \co(s))$, so it must be the complete linear series $|\co(s)|$; therefore, $\beta$ factors through the curve $\PP^1$ via the linear series $|\co(s)|$, i.e., $\Gamma$ lies on a rational normal curve: $\Gamma \to \PP^1 \xrightarrow{|\co(s)|} \PP^s$. The key constructions are summarized in the following diagram.

    \[\begin{tikzcd}[column sep=2em]
\PP^1 = \PP H^0(\PP^1, \co(1))^\vee &  & \PP H^0(\PP^1, \co(s))^\vee = \PP^s \\
      & \PP^1 \arrow[lu, "|\co(1)|", "\cong"'] \arrow[ru, hook, "|\co(s)|"'] &             \\
      & \Gamma \arrow[u] \arrow[luu, "\operatorname{res}(|\co(1)|) \subset |\co(1)|_\Gamma|", rounded corners=10pt, 
    to path={(\tikztostart) -| (\tikztotarget) \tikztonodes}] \arrow[ruu, "\operatorname{res}(W^\perp) \subset |\co(s)|_\Gamma|"', rounded corners=10pt, 
    to path={(\tikztostart) -| (\tikztotarget) \tikztonodes}] & 
\end{tikzcd}
\]
\end{proof}

Once we establish the existence of curves passing through prescribed points, a natural follow-up question is whether the embedding of the curve in the projective space is unique. As we discussed in Remark \ref{rmk: curves}, when the line bundles $L(-\Gamma)$ and $K_C \otimes L^\vee$ have negative degrees, Goppa duality induces a one-to-one relation on the maps to projective spaces. Moreover, because the factorization of the Gale dual pair through the curve $C$ by Goppa duality shares the same closed embedding $\Gamma \hookrightarrow C$, the problem reduces to whether the embedding of the curve in the projective space $\PP^s$ induced by the complete linear series $|\co(s)|$ fixing the mapping of the $s + 3$ points is unique. This follows from the fact that the only map fixing $s+2$ points on $\PP^1$ is the identity map. As a consequence, we obtain a proof for the classical result on the uniqueness of the rational normal curve in $\PP^s$ passing through $s + 3$ general points, summarized below.

.
\begin{proposition} \label{prop: uniqueness of RNC}
    Given $s + 3$ general points in $\PP^s$, there exists a unique rational normal curve passing through them.
\end{proposition}

\begin{proof}
    Let $\beta : \Gamma \to \PP^s$ be the closed embedding of $s + 3$ general points. Denote the Gale transform by $\alpha : \Gamma \to \PP^1$.  By Lemma \ref{lem: existence of RNC}, there exists a factorization of $\beta$ through a rational normal curve as $\Gamma \xrightarrow{\alpha} \PP^1 \xrightarrow{ \pi}\PP^s$.
    Suppose there exists another factorization of $\beta : \Gamma \to \PP^s$ through a rational normal curve, $\Gamma \xrightarrow{\alpha'} \PP^1 \xrightarrow{ \pi'}\PP^s$. By Goppa duality for curves applied to $\beta$, the Gale transform $\alpha$ of $\beta$ factors through $\alpha'$ as $\Gamma \xrightarrow{\alpha'} \PP^1 \xrightarrow{|\co(1)|} \PP^1$. It follows that $\alpha$ and $\alpha'$ differ by an automorphism of $\PP^1$, so by absorbing the automorphism in $\pi'$, we may assume  $\alpha' = \alpha$.  It remains to show $\pi'  = \pi$. 
    Because $\pi$ and $\pi'$ are induced by the same $(s+1)$-dimensional linear series, they differ by an automorphism of $\PP^s$; then the automorphism must fix $s +3$ general points on $\PP^s$, but the only such automorphism is the identity map. Hence $\pi' = \pi$, so the rational normal curve through $\Gamma$ is unique. 
\end{proof}

A similar technique can be applied to curves that are not the projective line, by using other uniqueness properties of curves through prescribed points. For example, if $r = 2$, then on one side of the duality, we are concerned with a map of the form $C \to \PP^2$ for some Gorenstein curve $C$. Recall hat there exists a unique degree-$d$ curve passing through $\frac{d(d+3)}{2}$ general points in $\PP^2$. Using this fact, if $\Gamma$ consists of $\gamma = \frac{d(d+3)}{2}$ points for some positive integer $d$ in $\PP^2$, then we can find a unique smooth curve $C$ passing through $\Gamma$, which fixes a factorization $\Gamma \to C \to \PP^2$. It follows that by Goppa duality, the Gale transform of $\Gamma \to \PP^2$ factors through this curve $C$ as $\Gamma \to C \to \PP^{\gamma-4}$. For example, suppose $\gamma = 9$ so that $d = 3$; then given nine general points in $\PP^2$, the unique cubic curve passing through them in $\PP^2$ is carried by Goppa duality to a curve in $\PP^5$, which passes through the Gale transform of the nine general points in $\PP^5$. Reversing the duality, we obtain the following result by Coble. 

\begin{proposition}[Coble {\cite[Theorem 19]{Coble1922}}]
    \label{prop: unique sectic elliptic cruve through 9 general points in P5}
    Given nine general points in $\PP^5$, there exists a unique sextic elliptic curve through them. 
\end{proposition}

\begin{proof}
    Let $\beta : \Gamma \to \PP^5$ be the map of the nine general points and denote its Gale transform by $\alpha : \Gamma \to \PP^2$. Since $\Gamma$ is general, there exists a unique elliptic curve $C$ in $\PP^2$ passing through the image of $\Gamma$, so $\alpha$ factors as 
    $\Gamma \xhookrightarrow{\iota} C \xhookrightarrow{\kappa} \PP^2,$
    where the second map $\kappa$ is a closed embedding induced by the complete linear series of a degree three line bundle $L$ on $C$. By Goppa duality for curves, $\beta$ factors through $C$ via a closed embedding $C \to \PP^5$ induced by the complete linear series of the degree six line bundle $K_C(\Gamma)\otimes L^\vee$.
    Therefore, there exists a sextic elliptic curve passing through the nine general points in $\PP^5$. 

    Now, suppose that there exists another factorization $\Gamma \xrightarrow{\iota'}
    C' \xrightarrow{\pi'} \PP^5$ of $\beta$ by a sextic elliptic curve $C'$. By a similar argument as above, $\pi'$ is induced by the complete linear series, and Goppa duality implies that $\alpha$ factors as $\Gamma \xrightarrow{\iota'}
    C' \xrightarrow{\kappa'} \PP^2$, where $\kappa$ is a closed embedding induced by the complete linear series of a degree three line bundle on $C'$. By uniqueness of an elliptic curve passing through nine general points, $\kappa$ and $\kappa'$ describe the same closed embedding, so we have an isomorphism of curves $C' \cong C$, and consequently  $\iota ' $ and $ \iota$ are the same under this isomorphism. It remains to show $\pi' = \pi$. Because $\pi'$ and $\pi$ are induced by the same complete linear series, they differ by an automorphism of $\PP^5$, which must fix the nine points $\Gamma$. Such an automorphism must be the identity map, hence $\pi' = \pi$ and the sextic elliptic curve is unique. They key constructions are summarized in the following diagram.

    
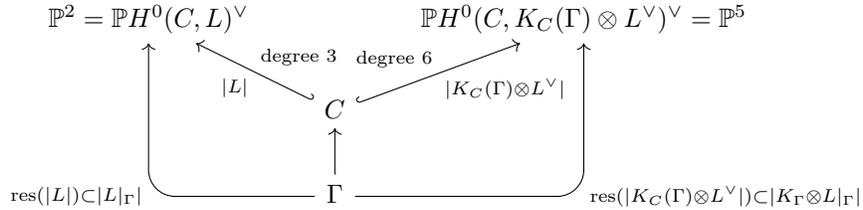
\[\begin{tikzcd}[column sep=2em]
\PP^2 = \PP H^0(C, L)^\vee &  & \PP H^0(C, K_C(\Gamma)\otimes L^\vee)^\vee = \PP^5 \\
      & C \arrow[lu, hook', "|L|", "\text{degree 3}"'] \arrow[ru, hook, "|K_C(\Gamma)\otimes L^\vee|"', "\text{degree 6}"] &             \\
      & \Gamma \arrow[u] \arrow[luu, "\operatorname{res}(|L|) \subset |L|_\Gamma|", rounded corners=10pt, 
    to path={(\tikztostart) -| (\tikztotarget) \tikztonodes}] \arrow[ruu, "\operatorname{res}(|K_C(\Gamma) \otimes L^\vee|) \subset |K_\Gamma \otimes L|_\Gamma|"', rounded corners=10pt, 
    to path={(\tikztostart) -| (\tikztotarget) \tikztonodes}] & 
\end{tikzcd}
\]
\end{proof}

\section{Goppa Duality For Surfaces}
\label{sec: 3}

Goppa duality for curves can be regarded as an extension of Gale duality by introducing a curve $C$ containing points $\Gamma$ as an effective Cartier divisor and promoting the line bundle $L$ to be defined on $C$. We wish to extend Goppa duality for curves (Theorem \ref{thm: Goppa duality for curves}) to higher dimensional schemes in a similar manner:  given a Gale dual pair, we want to factor the maps through a higher dimensional scheme $B$, whose maps to projective spaces are induced by linear series which are the annihilators of each other under Serre duality on $\Gamma$ after the restriction maps. To achieve this, we invoke a theorem proven by Eisenbud and Popescu \cite{EisenbudPopescu2000}, which gives us an exact sequence involving global sections of line bundles and the dual spaces, that resemble the derivation of Goppa duality for curves, as in the next subsection.

\subsection{Goppa Duality for Higher 
Dimensional Locally Gorenstein Schemes}
\label{sec: 3.1}

In the derivation of Goppa duality for curves, we used the fact that $\Gamma$ is an effective Cartier divisor on $C$, so that for any line bundle $L$ on $C$, we have a short exact sequence $0 \to L(-\Gamma) \to L \to L|_\Gamma \to 0$; we may regard this as a locally free resolution of the line bundle $L|_\Gamma$. Recall that the short exact sequence of sheaves induces a long exact sequence of cohomology, from which we determine the linear series $W$ of $L$ and $W^\perp$ of $K_C(\Gamma) \otimes L^\vee$ satisfying the short exact sequence that describes a Gale transform pair, as shown in the two-row diagram:
\bigskip

\begin{center}
\adjustbox{scale=0.9, center}{
\begin{tikzcd}[column sep=small]
{H^0(C, L(-\Gamma))} \arrow[r] & {H^0(C,L)} \arrow[r, "\operatorname{res}"] 
  & {H^0(\Gamma, L|_\Gamma) \cong H^0(\Gamma, K_\Gamma \otimes L|_\Gamma^\vee)^\vee} 
    \arrow[r, "\operatorname{res}^\vee"] \arrow[d, equal] 
  & {H^0(C,K_C(\Gamma) \otimes L^\vee)^\vee} \arrow[r] \arrow[d, two heads] 
  & {H^0(C,K_C \otimes L^\vee)^\vee} \\
0 \arrow[r] & W \arrow[u, hook] \arrow[r, "\operatorname{res}"] 
  & {H^0(\Gamma, L|_\Gamma) \cong H^0(\Gamma, K_\Gamma \otimes L|_\Gamma^\vee)^\vee} 
    \arrow[r, "\operatorname{res}^\vee"] 
  & (W^\perp)^\vee \arrow[r] & 0
\end{tikzcd}
\ \text{($\star$)}  
}
\end{center}
\bigskip

\noindent A higher-dimensional version of Goppa duality that we derive in Theorem \ref{thm: Goppa duality} similarly uses a resolution of a line bundle on $\Gamma$ by locally free sheaves on a higher dimensional scheme $B$. The following result by Eisenbud and Popescu shows the existence of an exact sequence of global sections and their dual spaces that have a similar structure as $(\star)$, which we will use to describe Goppa duality in more generality. 

\begin{theorem}[Eisenbud-Popescu \cite{EisenbudPopescu2000}, Theorem 3.1]
    \label{thm: 3.1}
    Let $\iota : \Gamma \to B$ be a finite map from a zero-dimensional Gorenstein scheme  to a Gorenstein scheme $B$ of dimension $c$, and let $L_\Gamma$ be a line bundle on $\Gamma$. Suppose that 
    $$0 \to \mathcal{E}_c \to  \mathcal{E}_{c-1} \to \cdots \to  \mathcal{E}_0 \to \iota_*L_\Gamma \to 0 $$
    is a resolution of $\iota_*L_\Gamma$  by locally free sheaves on $B$. Then $$\iota_*(K_\Gamma\otimes L_\Gamma^\vee) = \operatorname{coker} \big(\Hom(\mathcal{E}_{c-1} , K_B) \to  \Hom(\mathcal{E}_{c} , K_B) \big),$$ and if 
    $$H^{i+1} (B,\mathcal{E}_i) =H^{i+1} (B,\mathcal{E}_{i+1}) = 0 \text{ for all } 0 \leq i \leq c - 2, $$
    then the induced sequence

    \bigskip
    \adjustbox{scale=0.9, center}{%
\begin{tikzcd}[column sep=small]
{H^0(B, \mathcal{E}_1)} 
  \arrow[r] & {H^0(B, \mathcal{E}_0)} 
  \arrow[r, "\alpha"] & {H^0(\Gamma, L_\Gamma) \cong H^0(\Gamma, K_\Gamma \otimes L_\Gamma^\vee)^\vee} 
  \arrow[r, "\beta^\vee"] & {H^0(B, \Hom(\mathcal{E}_{c} , K_B))^\vee} 
  \arrow[r] & {H^0(B, \Hom(\mathcal{E}_{c-1} , K_B))^\vee}
\end{tikzcd}%
}
    \medskip
    \noindent is exact. 
\end{theorem}

Using the induced exact sequence, we derive Goppa duality for higher dimensional schemes as follows. We  find a pair $(W,W^\perp)$ of subspaces  of $H^0(B,\mathcal{E}_0)$ and $H^0(B, \Hom(\mathcal{E}_c, K_B))$ in a similar manner as in the derivation of Goppa duality for curves: $W$ is a subspace of $H^0(B,\mathcal{E}_0)$ complementary to the kernel of the map $ H^0(B, \mathcal{E}_0) \xrightarrow{\alpha} H^0(\Gamma, L_\Gamma)$, and similarly, $W^\perp$ is a subspace of $H^0(B,\Hom(\mathcal{E}_c, K_B))$ complementary to the kernel of the map $ H^0(B, \Hom(\mathcal{E}_c,K_B)) \xrightarrow{\beta}  H^0(\Gamma, L_\Gamma)$. It follows that we have a short exact sequence of vector spaces involving $W$ and $W^\perp$, which induces the Gale dual pair, generalizing $(\star)$:
\bigskip

\adjustbox{scale=0.9, center}{
\begin{tikzcd}[column sep=tiny]
{H^0(B, \mathcal{E}_1)} \arrow[r] & {H^0(B,\mathcal{E}_0)} \arrow[r, "\alpha"]        & {H^0(\Gamma, L_\Gamma) \cong H^0(\Gamma, K_\Gamma \otimes L_\Gamma^\vee)^\vee} \arrow[r, "\beta^\vee"] \arrow[d, equal] & {H^0(B, \Hom(\mathcal{E}_c, K_B))^\vee} \arrow[r] \arrow[d, two heads] & {H^0(B,\Hom(\mathcal{E}_{c-1}, K_B))^\vee} \\
0 \arrow[r]                    & W \arrow[u, hook] \arrow[r, "\operatorname{}"] & {H^0(\Gamma, L_\Gamma) \cong H^0(\Gamma, K_\Gamma \otimes L_\Gamma^\vee)^\vee} \arrow[r, "\operatorname{}"]                                & (W^\perp)^\vee \arrow[r]                                                       & 0   .
\end{tikzcd}
}
\bigskip

\noindent In the case $\mathcal{E}_0$ and $\Hom(\mathcal{E}_c, K_B)$ are line bundles on the base scheme $B$, we have the following algebro-geometric interpretation of the Gale dual pair: the maps $\Gamma \to \PP W^\vee$ and $\Gamma \to \PP (W^\perp)^\vee$  factor through $B$ via the linear series $W$ and $W^\perp$, respectively. We call this result \emph{generalized Goppa duality}, or simply \emph{Goppa duality}, as classical Goppa duality introduced in Theorem \ref{thm: Goppa duality for curves} is a special version of this result. 

\begin{theorem} [generalized Goppa duality]
    \label{thm: Goppa duality}
    Let $\iota: \Gamma \to B$ be a closed embedding of a finite Gorenstein scheme $\Gamma$ of length $\gamma$ into a Gorenstein scheme $B$, of codimension $c$. Suppose that $L_\Gamma$ is a line bundle on $\Gamma$ with a resolution 
    $$0 \to \mathcal{E}_c \to \mathcal{E}_{c-1} \to \cdots \to \mathcal{E}_0 \to \iota_*L_\Gamma \to 0$$
    by locally free sheaves on $B$, where $\mathcal{E}_0$ and $\mathcal{E}_c$ are line bundles.
    Assume that the sheaves satisfy
    $$H^{i+1} (B,\mathcal{E}_i) =H^{i+1} (B,\mathcal{E}_{i+1}) = 0 \text{ for all } 0 \leq i \leq c - 2. $$
    
    \noindent Suppose that $W$ is an $(r+1)$-dimensional linear series of $\mathcal{E}_0$ complementary to the kernel of the map $H^0(B, \mathcal{E}_0) \xrightarrow{\alpha} H^0(\Gamma, L_\Gamma)$. Then there exists an $(s+1)$-dimensional linear series $W^\perp$ of $K_B \otimes \mathcal{E}_c^\vee$, satisfying the following short exact sequence,
    $$0 \to W \xrightarrow{} H^0(\Gamma, L_\Gamma) \cong H^0(\Gamma, K_\Gamma\otimes L_\Gamma^\vee)^\vee \xrightarrow{} (W^\perp)^\vee \to 0,$$
    and such a linear series $W^\perp$ is complementary to the kernel of the map $H^0(B, K_B \otimes \mathcal{E}_c^\vee) \xrightarrow{\beta} H^0(\Gamma, K_\Gamma\otimes L_\Gamma^\vee)$.
    
    In particular, the induced Gale dual pair $\Gamma \to \PP^r$ and $\Gamma \to \PP^s$ factors through $B$ as 
    $$\Gamma \xrightarrow{\iota} B \xdashrightarrow{W} \PP^r \text{ and } \Gamma \xrightarrow{\iota} B \xdashrightarrow{W^\perp} \PP^s, $$
    where the rational maps $B \xdashrightarrow{W} \PP^r$ and $ B \xdashrightarrow{W^\perp} \PP^s$ are induced by the linear series $W$ and $W^\perp$, respectively.
    
\end{theorem}

\begin{proof}
By hypothesis of vanishing of the first cohomology group $H^1$ and Theorem \ref{thm: 3.1}, we obtain the following exact sequences:
\medskip

\adjustbox{scale=0.9, center}{
\begin{tikzcd}[column sep=tiny]
{H^0(B, \mathcal{E}_1)} \arrow[r] & {H^0(B,\mathcal{E}_0)} \arrow[r, "\alpha"]        & {H^0(\Gamma, L_\Gamma) \cong H^0(\Gamma, K_\Gamma \otimes L_\Gamma^\vee)^\vee} \arrow[r, "\beta^\vee"] \arrow[d, equal] & {H^0(B,  K_B \otimes \mathcal{E}_c^\vee)^\vee} \arrow[r] \arrow[d, two heads] & {H^0(B, K_B \otimes \mathcal{E}_{c-1}^\vee)^\vee} \\
0 \arrow[r]                    & W \arrow[u, hook] \arrow[r] & {H^0(\Gamma, L_\Gamma) \cong H^0(\Gamma, K_\Gamma \otimes L_\Gamma^\vee)^\vee} \arrow[r, ""]                                & (W^\perp)^\vee \arrow[r]                                                       & 0 .                              
\end{tikzcd}
}
\medskip

\noindent It follows that the induced Gale dual pair $\Gamma \to \PP^r$ and $\Gamma \to \PP^s$ factors through $B$ via the linear series $W$ and $W^\perp$, respectively. The key constructions are summarized in the following diagram.

\[\begin{tikzcd}[column sep=6em]
\PP^r = \PP W^\vee &  & \PP (W^\perp)^\vee = \PP^s \\
      & B \arrow[lu, dashed, "W \subset |L|"] \arrow[ru, dashed, "W^\perp \subset |K_B \otimes \mathcal{E}_c^\vee|"'] &             \\
      & \Gamma \arrow[u] \arrow[luu, "\alpha(W) \subset |L|_\Gamma|", rounded corners=10pt, 
    to path={(\tikztostart) -| (\tikztotarget) \tikztonodes}] \arrow[ruu, "\beta(W^\perp) \subset |K_\Gamma \otimes L_\Gamma^\vee|"', rounded corners=10pt, 
    to path={(\tikztostart) -| (\tikztotarget) \tikztonodes}] & 
\end{tikzcd}
\]
\end{proof}

\subsection{When $\Gamma$ is a Complete Intersection on $S$}
\label{sec: complete intersection}

For the rest of the paper, we apply generalized Goppa duality to the case in which the base scheme $B$ is a surface $S$. Because Theorem \ref{thm: Goppa duality} requires a resolution of a line bundle on $\Gamma$, we focus on the case in which such a resolution exists canonically: when $\Gamma$ is a complete intersection on the surface $S$, and $L_\Gamma$ is the restriction of a line bundle $L$ on $S$. In this case, we can use the Koszul resolution of the structure sheaf of the complete intersection to induce a resolution of $L|_\Gamma$ for any line bundle $L$ on $S$, and Goppa duality applies.
We begin by recalling the Koszul resolution for a complete intersection of two divisors on a surface.

\begin{lemma}[Eisenbud \cite{Eisenbud}]
    \label{lem: loc free resolution of complete intersection}
    Let $\Gamma$ be a complete intersection of two effective Cartier divisors $C_1, C_2$ on a  surface $S$. Then the ideal sheaf $\mathcal{I}$ of $\Gamma$ admits a locally free Koszul resolution
    $$0 \to \co_S(-C_1 - C_2) \to \co_S(-C_1) \oplus \co_S(-C_2) \to \mathcal{I} \to 0. $$ 
\end{lemma} 
\noindent After twisting this exact sequence by $L$, the locally free resolution of the ideal sheaf of $\Gamma$ in $S$ induces a length three locally free resolution of $L|_\Gamma$:
$$0 \to L(-C_1 -C_2) \to L(-C_1) \oplus L(-C_2) \to L \to L|_\Gamma \to 0.$$
Goppa duality (Theorem \ref{thm: Goppa duality}) can then be applied, where $\alpha$ and $\beta$ correspond to restriction maps of sections of the corresponding line bundles on $C$ to $\Gamma$. This leads to the following result of Goppa duality for a complete intersection on a surface.

\begin{theorem}[Goppa Duality for a Complete Intersection on a Surface]
    \label{thm: goppa duality for surfaces}
    Let $\iota: \Gamma \to S$ be the closed embedding of a length-$\gamma$ complete intersection of two effective Cartier divisor curves $C_1, C_2$ on a Gorenstein surface $S$. Assume that $L$ is a line bundle on $S$ such that 
    $$H^1(S, L) = H^1(S, L(-C_1)) = H^1(S, L(-C_2)) = 0.$$
    Suppose that $W$ is an $(r+1)$-dimensional linear series  of $L$  complementary to the kernel of the restriction map $H^0(C,L) \xrightarrow{\operatorname{res}}H^0(C,L|_\Gamma)$. Then there exists an $(s+1)$-dimensional linear series $ W^\perp$ of $K_S(C_1 + C_2)\otimes L^\vee$ satisfying the following short exact sequence,
    $$0 \to W \xrightarrow{\operatorname{res}} H^0(\Gamma, L|_\Gamma) \cong H^0(\Gamma, K_\Gamma\otimes L|_\Gamma^\vee)^\vee \xrightarrow{\operatorname{res}^\vee} (W^\perp)^\vee \to 0,$$
    and such a linear series $W^\perp$ is complementary to the kernel of the restriction map $H^0(S, K_S(C_1+C_2) \otimes L^\vee)\xrightarrow{\operatorname{res}} H^0(\Gamma, K_\Gamma\otimes L|_\Gamma^\vee)$.
    
    In particular, the Gale dual pair $\Gamma \to \PP^r$ and $\Gamma \to \PP^s$ factors through $S$ as 
    $$\Gamma \xrightarrow{\iota} S \xdashrightarrow{W} \PP^r \text{ and } \Gamma \xrightarrow{\iota} S \xdashrightarrow{W^\perp} \PP^s. $$
\end{theorem}

\noindent We will refer to this result as \textit{Goppa duality for a complete intersection on a surface} or \emph{Goppa duality for surfaces} in short. We will also call the line bundle $K_S(C_1 + C_2) \otimes L^\vee$ the \textit{Goppa dual line bundle} of $L$, and the linear series  $W^\perp$ the \textit{Goppa dual linear series} of $W$.

As a restatement, Goppa duality for surfaces is a duality between two maps $$S \xdashrightarrow{W} \PP^r \text{ and }  S \xdashrightarrow{W^\perp} \PP^s,$$
where $W$ is a linear series of $L$ and $W^\perp$ is the linear series of $K_S(C_1+C_2) \otimes L^\vee$, whose restrictions to $\Gamma$ induce an orthogonal pair under Serre duality on $\Gamma$. We may understand Goppa duality in relation to Gale duality by filling in the leveled diagram as below:
\[\begin{tikzcd}[column sep=6em, row sep=small]
\PP W^\vee &  & \PP (W^\perp)^\vee \\
      & S \arrow[lu, dashed, "W \subset |L|"] \arrow[ru, dashed, "W^\perp \subset |K_S(C_1+C_2) \otimes L^\vee|"'] &             \\
      & C_i \arrow[u] &             \\
      & \Gamma \arrow[u] \arrow[luuu, "\operatorname{res}(W) \subset |L|_\Gamma|", rounded corners=10pt, 
    to path={(\tikztostart) -| (\tikztotarget) \tikztonodes}] \arrow[ruuu, "\operatorname{res}(W^\perp) \subset |K_\Gamma \otimes L|_\Gamma^\vee|"', rounded corners=10pt, 
    to path={(\tikztostart) -| (\tikztotarget) \tikztonodes}] & 
\end{tikzcd}
\]

\bigskip

\begin{remark}
\label{ref: generalized goppa duality}
In parallel to Remark \ref{rmk: curves}, we emphasize that the choice of the linear series $W$ or $W^\perp$ is not unique in general. Nevertheless, when both of the restriction maps of global sections of $L$ and $K_S(C_1+C_2) \otimes L^\vee$ on $C$ to $\Gamma$ are injective,  Goppa duality for surfaces induces a one-to-one correspondence between the two maps  $S \xdashrightarrow{|L|} \PP^r$ and  $S \xdashrightarrow{|K_S(C_1 +C_2) \otimes L^\vee|} \PP^s$ factoring the Gale dual pair.

\end{remark}

Our goal for the rest of this paper is to use Goppa duality for surfaces to deduce existence or uniqueness statements of surfaces passing through  prescribed points in projective spaces, in a similar manner as in Section \ref{sec: Goppa duality}. Recall that in the discussion of existence and uniqueness of rational normal curves, we chose the curve $C$ to be the projective line and set $r = 1$, so that the morphism $C \to \PP^1$ in the factorization on one side of the Goppa duality can be taken to be the identity map. A similar approach can be applied to Goppa duality for surfaces, by taking $S$ to be the projective plane and $r = 2$, and regarding the map $S \to \PP^2$ to be the identity map, which is to say that $d = 1$ and $W$ is the complete linear series $|\co(1)|$. We develop this discussion further in the following section.

\subsection{Case $S = \PP^2$}
\label{sec: S=PP2}

In this section, we focus on the case in which $S$ is the projective plane $\PP^2$. This setting makes Goppa duality particular useful: the vanishing requirement of the first cohomology group $H^1$ in Theorem \ref{thm: 3.1} is automatically satisfied for any line bundle.  Moreover, line bundles on $\PP^2$ are of the form $\co(d)$ for some integer $d$, so the computation of the Goppa dual line bundle $K_S(C_1+C_2) \otimes L^\vee$ is straightforward.  For example, the Goppa dual linear series arising from the Goppa duality for a $(d_1, d_2)$-complete intersection on $\PP^2$ embedded via $\mathcal{O}(d)$ corresponds to a linear series of  $\co(d_1 + d_2 - d - 3)$. This leads to the following version of Goppa duality for surfaces, which we will call \emph{Goppa duality for a complete intersection on $\PP^2$} or simply \emph{Goppa duality for $\PP^2$}.

\begin{corollary}[Goppa Duality for a Complete Intersection on $\PP^2$]
    \label{cor:complete intersection on P2}
    Let $\iota: \Gamma \to \PP^2$ be the closed embedding of a $(d_1,d_2)$-complete intersection on $\PP^2$. 
    Suppose that $W$ is an $(r+1)$- dimensional linear series of the line bundle $\co(d)$ complementary to the kernel of the restriction map $H^0(\PP^2, \co(d)) \to H^0(\Gamma, \co(d)|_\Gamma)$. Then there exists an $(s+1)$-dimensional linear series $W^\perp$ of  $ \co(d_1+d_2-d-3)$ satisfying the following short exact sequence:
    $$0 \to W \xrightarrow{\operatorname{res}} H^0(\Gamma, L|_\Gamma) \cong H^0(\Gamma, K_\Gamma\otimes L|_\Gamma^\vee)^\vee \xrightarrow{\operatorname{res}^\vee} (W^\perp)^\vee \to 0,$$
    and such a linear series $W^\perp$ is complementary to the kernel of the restriction map $H^0(\PP^2, \co(d_1+d_2-d-3)) \to H^0(\Gamma, \co(d_1+d_2-d-3)|_\Gamma)$.  
    
    In particular, the Gale dual pair $\Gamma \to \PP^r$ and $\Gamma \to \PP^s$ factors through $\PP^2$ as 
    $$\Gamma \xrightarrow{\iota} \PP^2 \xdashrightarrow{W} \PP^r \text{ and } \Gamma \xrightarrow{\iota} \PP^2 \xdashrightarrow{W^\perp} \PP^s. $$
    In terms of our level diagram relating Goppa duality to Gale duality, we have the following diagram.
    \[\begin{tikzcd}[column sep=3em, row sep=small]
\PP^r = \PP W^\vee &  & \PP (W^\perp)^\vee = \PP^s \\
      & \PP^2 \arrow[lu, dashed, "W \subset |\co(d)|"] \arrow[ru, dashed, "W^\perp \subset |\co(d_1+d_2-d-3)|"'] &             \\
      & C_i \arrow[u] &             \\
      & \Gamma \arrow[u] \arrow[luuu, "\operatorname{res}(W) \subset |\co(d)|_\Gamma|", rounded corners=10pt, 
    to path={(\tikztostart) -| (\tikztotarget) \tikztonodes}] \arrow[ruuu, "\operatorname{res}(W^\perp) \subset |\co(d_1+d_2-d-3)|_\Gamma|"', rounded corners=10pt, 
    to path={(\tikztostart) -| (\tikztotarget) \tikztonodes}] & 
\end{tikzcd}
\]

\bigskip

\end{corollary}

In a similar manner as in Section \ref{sec: Goppa duality}, in practice, we use Goppa duality for $\PP^2$ as follows: given $\gamma$ points $\Gamma \to \PP^r$ in a projective space, we try to realize $\Gamma$ as a complete intersection on $\PP^2$ satisfying the following conditions:
\begin{enumerate}[(i)]
    \item $\PP^2$ factors the map $\Gamma \to \PP^r$ via a closed embedding $\iota: \Gamma \hookrightarrow \PP^2$, through which $\Gamma$ can be regarded as a $(d_1, d_2)$-complete intersection, and
    \item the map $\PP^2 \dashrightarrow \PP^r$ in the factorization is induced by an $(r+1)$-dimensional linear series $W$ of $\co_{\PP^2}(d)$ complementary to with the kernel of the restriction map $H^0(\PP^2, \co(d)) \to H^0(\Gamma, \co(d)|_\Gamma)$. 
\end{enumerate} 

\noindent Once we find a factoring of the original map through $\PP^2$ satisfying these conditions, we can factor the Gale transform $\Gamma \to \PP^s$ via $\PP^2$ as the composition of the closed embedding $\iota$ and the map $C \dashrightarrow \PP^s$, which is induced by a $(s+1)$-dimensional linear series $W^\perp$.


\subsubsection{First example: $(3,3)$-complete intersection on $\PP^2$}
\label{sec: (3,3)-complete intersection}

As the first nontrivial example of Goppa duality for surfaces applied to the case $S = \PP^2$ and $r = 2$, consider the case where on the other side of the duality, the map $\PP^2 \dashrightarrow \PP^s$ from the surface to projective space is the Veronese embedding $\PP^2 \hookrightarrow \PP^5$. This implies that $s = 5$ and the Goppa dual linear series $W^\perp$ is required to be the complete linear series $|\co(2)|$. For Goppa duality for $\PP^2$ to hold, $\Gamma$ needs to be a $(d_1, d_2)$-complete intersection of length $r+s + 2 =9$, and $\co(d_1 + d_2 - 1 - 3)= \co(2)$. The only possible situation satisfying these conditions is $d_1=d_2 = 3,$ i.e., when $\Gamma$ is a $(3,3)$-complete intersection on $\PP^2$. Conversely, because the restriction maps $H^0(\PP^2, \co(1)) \to  H^0(\Gamma, \co(1)|_\Gamma)$ and $H^0(\PP^2, \co(2)) \to  H^0(\Gamma, \co(2)|_\Gamma)$ are both injective for a $(3,3)$-complete intersection, the linear series satisfy the requirement of Goppa duality for surfaces. In summary, if $\Gamma \to \PP^2$ is a general $(3,3)$-complete intersection, then its Gale transform in $\PP^5$ factors through a Veronese surface.

\begin{lemma}
    \label{lem: (3,3)-CI existence}
    Suppose that $\Gamma \xrightarrow{\alpha} \PP^2$  is the closed embedding of a general $(3,3)$-complete intersection. Then its Gale transform factors through a Veronese surface: $\Gamma \xrightarrow{\alpha} \PP^2 \xhookrightarrow{|\mathcal{O}(2)|} \PP^5$.
\end{lemma}

\begin{proof}
    By Goppa duality for $\PP^2$ (Corollary \ref{cor:complete intersection on P2}) with the choice $(d,d_1,d_2) = (1,3,3)$ and $W = |\co(1)|$, conditions (i) and (ii) are satisfied, so the Gale transform $\beta: \Gamma \to \PP^5$ of $\alpha$ factors through the surface $\PP^2$ as $\Gamma \xrightarrow{\alpha} \PP^2 \xrightarrow{\pi} \PP^5$, where $\pi$ is induced by a 6-dimensional linear series $W^\perp$ of the Goppa dual line bundle $\co(3+3-1-3) = \co(2)$. In particular, $W^\perp$ is the complete linear series, so the map $\PP^2 \to \PP^5$ factoring the Gale transform $\Gamma \to \PP^5$ is a Veronese embedding.
\end{proof}


Furthermore, an analogous technique as in the proof of uniqueness of the rational normal curve in Proposition \ref{prop: uniqueness of RNC} can be applied to this setting, so that the factorization through the Veronese surface in this example is unique: we compare two factorizations through a Veronese surface, and invoke rigidity of projective automorphisms that fix a sufficiently large number of general points. The full result and proof are summarized in the following Lemma.

\begin{lemma}
    \label{lem: (3,3)-CI uniqueness}
    Suppose that $\Gamma \xrightarrow{\alpha} \PP^2$  is the closed embedding of a general $(3,3)$-complete intersection. Then the Gale transform factors through a Veronese surface $\PP^2$ as $\Gamma \xrightarrow{\alpha} \PP^2 \xrightarrow{|\mathcal{O}(2)|} \PP^5$, and the last map is unique up to reparametrization of $\PP^2$.
\end{lemma}

\begin{proof}
    Suppose that we have another factorization of the Gale transform through a Veronese surface as $\Gamma \xrightarrow{\alpha'} \PP^2 \xrightarrow{ \pi'} \PP^5$, where $\alpha'$ describes the closed embedding of another $(3,3)$-complete intersection. Then the map $\alpha$ factors through $\PP^2$ as $\Gamma \xrightarrow{\alpha'} \PP^2 \xrightarrow{|\mathcal{O}(1)| , \kappa} \PP^2$ by Goppa duality $\PP^2$ (Corollary \ref{cor:complete intersection on P2}). Because $\kappa$ is an automorphism of $\PP^2$, we may reduce to the case where $\alpha' = \alpha$. To show $\pi' = \pi$, observe that $\pi$ and $\pi'$ differ by an automorphism of $\PP^5$ since they are induced by the same linear series. Because an automorphism of $\PP^5$ can fix at most seven general points in $\PP^5$ and $\Gamma$ is a general $(3,3)$-complete intersection,  $\phi$ must be the identity. Therefore, the factorization is unique up to a reparametrization of $\PP^2$. 
    Key aspects of the argument are summarized in the following diagram.
    
    \[\begin{tikzcd}[column sep=2em, row sep=small]
\PP^2 = \PP H^0(\PP^2,\co(1))^\vee &  & \PP H^0(\PP^2,\co(2))^\vee = \PP^5 \\
      & \PP^2 \arrow[lu, "|\co(1)|", "\cong"'] \arrow[ru, hook, "W^\perp \subset |\co(2)|"'] &             \\
      & C_i \arrow[u] &             \\
      & \Gamma \arrow[u] \arrow[luuu, "\operatorname{res}(|\co(1)|) \subset |\co(1)|_\Gamma|", rounded corners=10pt, 
    to path={(\tikztostart) -| (\tikztotarget) \tikztonodes}] \arrow[ruuu, "\operatorname{res}(|\co(2)|) \subset |\co(2)|_\Gamma|"', rounded corners=10pt, 
    to path={(\tikztostart) -| (\tikztotarget) \tikztonodes}] & 
\end{tikzcd}
\]
\end{proof}

Note that this lemma assumes that the nine points are given in $\PP^2$.
In parallel to the argument of uniqueness of rational normal curves in $\PP^s$ through $s+3$ general points, we may hope to start with nine points in $\PP^5$ instead and deduce existence or uniqueness of a Veronese surface passing through the points. However, there is a caveat to this idea: if we are given only the map $\Gamma \to \PP^5$ of nine points in $\PP^5$, then the Gale transform $\Gamma \to \PP^2$ is not necessarily a $(3,3)$-complete intersection; we want to require that the points are special, lying on two cubics in $\PP^2$. Hence, to guarantee that the Gale transform in $\PP^2$ is a complete intersection, we need to prescribe curves in $\PP^5$ which would realize the Gale transform $\Gamma \to \PP^2$ as a complete intersection on $\PP^2$ via Goppa duality for curves. 
More precisely, suppose that we want to impose that on the $\PP^2$ of the Gale duality, the nine nondegenerate points lie on a $(3,3)$-complete intersection, say on two degree three curves $C_1$ and $C_2$. If we denote the line bundle $L_i$ on $C_i$ associated to the closed embedding $C_i \hookrightarrow \PP^2$, then Goppa duality for curves (Theorem \ref{thm: Goppa duality for curves}) shows that the Gale transform $\Gamma \to \PP^5$ factors through a map $C \to \PP^5$ induced by a nondegenerate linear series of the degree six line bundle $K_C(\Gamma) \otimes L_i^\vee$. By reversing the duality, it follows that if we start with two prescribed sextic elliptic curves in $\PP^5$ passing through the nine points in $\PP^5$, then  the Gale transform $\Gamma \to \PP^2$ factors through $C_i \to \PP^2$ of degree three for each $i$, which implies that $\Gamma$ is a $(3,3)$-complete intersection on $\PP^2$, as desired. The key constructions are summarized in the following diagram.

    \[\begin{tikzcd}[column sep=3em]
\PP^2 = \PP H^0(C_i, L_i)^\vee &  & \PP H^0(C_i, K_{C_i}(\Gamma)\otimes L_i^\vee)^\vee = \PP^5 \\
      & C_i \arrow[lu, hook', dashed, "\text{induced}", "\text{degree 3}"'] \arrow[ru, hook, "\text{prescribed}"', "\text{degree 6}"] &             \\
      & \Gamma \arrow[u] \arrow[luu, "", rounded corners=10pt, 
    to path={(\tikztostart) -| (\tikztotarget) \tikztonodes}] \arrow[ruu, ""', rounded corners=10pt, 
    to path={(\tikztostart) -| (\tikztotarget) \tikztonodes}] & 
\end{tikzcd}
\]

\noindent In summary, once two sextic elliptic curves through $\Gamma$ are prescribed in $\PP^5$, Goppa duality for $\PP^2$  applies and we can induce a unique Veronese surface containing $\Gamma$ in $\PP^5$.

\begin{proposition}
    \label{prop: unique Veronese through two sextic elliptic curves}
    Suppose that $\Gamma$ consists of nine nondegenerate points in $\PP^5$ and there exist exactly two independent nondegenerate sextic elliptic curves $C_1, C_2$ passing through $\Gamma$ in $\PP^5$. Then there exists a unique Veronese surface containing $C_1$ and $C_2$.
\end{proposition}

\begin{proof}
    By the nondegeneracy assumption, the closed embedding $C_i \hookrightarrow \PP^5$ of each sextic curve is induced by the complete linear series of a degree 6 line bundle, denoted by $L_i$. By Goppa duality for curves (Theorem \ref{thm: Goppa duality for curves}), the Gale transform $\Gamma \to \PP^2$ factors through each $C_i$ as $\Gamma \to C_i \xrightarrow{K_{C_i}(\Gamma)\otimes L_i^\vee} \PP^2$, where $\deg(K_{C_i}(\Gamma)\otimes L_i^\vee) = 3$, and consequently the last map is a nondegenerate closed embedding. It follows that the Gale transform $\Gamma \to \PP^2$ is a $(3,3)$-intersection in $\PP^2$, so by Goppa duality for $\PP^2$ (Corollary \ref{cor:complete intersection on P2}), the original map factors as $\Gamma \hookrightarrow C_i \hookrightarrow \PP^2 \xrightarrow{|\mathcal{O}(2)|}\PP^5$, and this Veronese surface is unique up to reparametrization of $\PP^5$ by Lemma \ref{lem: (3,3)-CI uniqueness}. 
    The key relations of maps to projective spaces from the points, the curves, and the surfaces are summarized in the following level diagram.
    
    
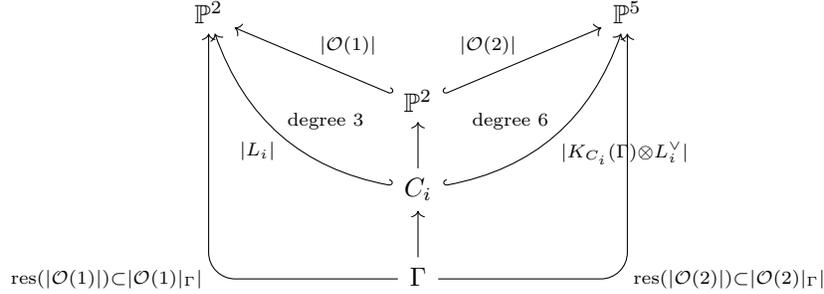
\[\begin{tikzcd}[column sep=6em]
\PP^2 &  & \PP^5 \\
      & \PP^2 \arrow[lu, hook', "|\co(1)|"'] \arrow[ru, hook, "|\co(2)|"] &             \\
      & C_i \arrow[u] \arrow[luu, hook', "|L_i|", "\text{degree 3}"', bend left] \arrow[ruu, hook, "|K_{C_i}(\Gamma)\otimes L_i^\vee|"',  "\text{degree 6}", bend right] &             \\
      & \Gamma \arrow[u] \arrow[luuu, "\operatorname{res}(|\co(1)|) \subset |\co(1)|_\Gamma|", rounded corners=10pt, 
    to path={(\tikztostart) -| (\tikztotarget) \tikztonodes}] \arrow[ruuu, "\operatorname{res}(|\co(2)|) \subset |\co(2)|_\Gamma|"', rounded corners=10pt, 
    to path={(\tikztostart) -| (\tikztotarget) \tikztonodes}] & 
\end{tikzcd}
\]
\end{proof}

\noindent Forgetting the configuration of the nine points leads to the following geometric corollary on a surface containing prescribed curves.
\begin{corollary}
    \label{cor: unique Veronese through two deg 6 elliptic curves in P5}
    Two nondegenerate sextic elliptic curves $C_i$ in $\PP^5$ intersecting at exactly nine nondegenerate points lie on a unique Veronese surface.  $\qed$
\end{corollary}

\subsubsection{More general example: $(2,d)$- and $(3,d)$-complete intersections on $\PP^2$}

More generally, we can apply the approach used in Proposition \ref{prop: unique Veronese through two sextic elliptic curves} to general $(2,d)$-complete intersections in $\PP^2$ with $d \geq 3$ and to general $(3,d)$-complete intersections in $\PP^2$ with $d \geq 4$. In a similar manner as in Section \ref{sec: (3,3)-complete intersection}, it is necessary to prescribe special curves in the larger projective space, as described below. 

Consider the case in which on one side of Goppa duality, $\Gamma$ on $\PP^2$ is achieved as a  $(2,d)$-complete intersection.  Note that the condition $d \geq 3$ guarantees that the Gale transform in $\PP^2$ always has a degree $d$ curve passing through them. Hence we need only to prescribe a rational normal curve of degree $2(d-2)$ in $\PP^{2d-4}$ passing through the points; then using Goppa duality for curves, we can induce a special conic in $\PP^2$ that passes through the $2d$ points $\Gamma$. The key idea for inducing the special conic is shown the following diagram.

    
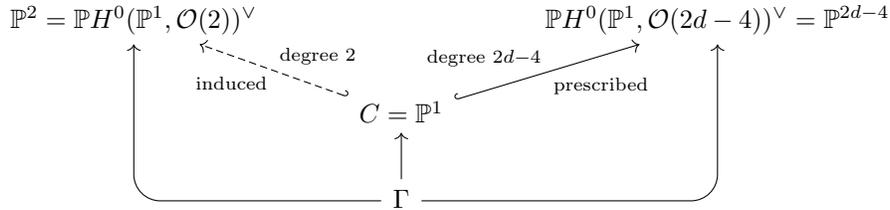
\[\begin{tikzcd}[column sep=3em]
\PP^2 = \PP H^0(\PP^1, \co(2))^\vee &  & \PP H^0(\PP^1, \co(2d-4))^\vee = \PP^{2d-4} \\
      & C = \PP^1 \arrow[lu, hook', dashed, "\text{induced}", "\text{degree 2}"'] \arrow[ru, hook, "\text{prescribed}"', "\text{degree }2d-4"] &             \\
      & \Gamma \arrow[u] \arrow[luu, "", rounded corners=10pt, 
    to path={(\tikztostart) -| (\tikztotarget) \tikztonodes}] \arrow[ruu, ""', rounded corners=10pt, 
    to path={(\tikztostart) -| (\tikztotarget) \tikztonodes}] & 
\end{tikzcd}
\]

\noindent Once the special curve is given, we apply Goppa duality for  $\PP^2$ (Corollary \ref{cor:complete intersection on P2}) to find a factorization of $\Gamma \to \PP^{2d-4}$ through the surface $\PP^2$. On the $\PP^2$ side of Goppa duality, we must set $W$ to be the complete linear series $|\co_{\PP^2}(1)|$ to satisfy condition (ii). The Goppa dual line bundle is then $\co_{\PP^2}(2 + d - 1 - 3) = \co_{\PP^2}(d-2)$, whose space of global sections have dimension $\frac{d(d-1)}{2}$. On the other hand, the Goppa dual linear series $W^\perp$ must be of dimension $2d-3$, being complementary to the kernel of the restriction map $H^0(\PP^2, \co(d-2)) \to H^0(\Gamma, \co(d-2)|_\Gamma)$. It follows that  we have a choice for the linear series $W^\perp$ satisfying condition (ii). 
To determine conditions for $W^\perp$ to satisfy this property, observe that the restriction map fits into the induced exact sequence of cohomology,
\begin{align*}
    0 &\to H^0(\PP^2, \mathcal{I}(d-2)) \to H^0(\PP^2, \co(d-2)) \xrightarrow{\operatorname{res}} H^0(\Gamma, \co(d-2)|_\Gamma) .
\end{align*}
Moreover, $H^0(\PP^2, \mathcal{I}(d-2))$ can be computed via a twist of the Koszul resolution, $$ 0 \to \co(-2-d) \to \co(-2)\oplus \co(-d) \to \mathcal{I} \to 0,$$
which induces the long exact sequence of cohomology
\begin{align*}
    \cancel{H^0(\PP^2, \co(-4))} \to H^0(\PP^2, \co(d-4) \oplus \co(-2)) \to H^0(\PP^2, \mathcal{I}(d-2)) 
    \to \cancel{H^1(\PP^2, \co(-4)) }.
\end{align*}
In particular, if $d \geq 4$, then we have an isomorphism $$H^0(\PP^2, \co(d-4)) \xrightarrow{\cong} H^0(\PP^2, \mathcal{I}(d-2)).$$ Geometrically, we can interpret this isomorphism as follows. Let $f$ be the degree two homogeneous equation of the special conic passing through $\Gamma$ in $\PP^2$. Under the restriction map, the homogeneous equation vanishes entirely on $\Gamma$ by construction. The isomorphism of global sections is then given by multiplication of degree -$(d-4)$ homogeneous polynomials by $f$. It follows that the kernel of the restriction map is isomorphic to $H^0(\PP^2, \co(d-4))\cdot f$. Therefore, the Goppa dual linear series $W^\perp$ must be chosen among the complementary subspaces to $H^0(\PP^2, \co(d-4))\cdot f$ in $H^0(\PP^2, \co(d-2))$, and converse also holds.   It follows that  the moduli space of $W^\perp$ is $\dim G(2d-3, d(d-1)/2)$-dimensional: the family of $(2d-3)$-dimensional vector subspace of $H^0(\PP^2, \co(d-2))$ is the Grassmannian $G(2d-3, d(d-1)/2)$, and the condition to be complementary to the vector subspace $H^0(\PP^2, \co(d-4))\cdot f$ is a proper closed condition in the Grassmannian.  In summary, we have the following level diagram describing the duality, and we obtain the following result on the existence of a factorization of a rational normal curve through a rational map $\PP^2 \dashrightarrow \PP^s$.

    \[\begin{tikzcd}[column sep=6em]
\PP^2 &  & \PP^{2d-4} \\
      & \PP^2 \arrow[lu, hook', "|\co_{\PP^2}(1)|", "\cong"'] \arrow[ru, dashed, "W \subset |\co_{\PP^2}(d-2)|"' , "\text{choices}"] &             \\
      & \PP^1 \arrow[u, hook, "|\co_{\PP^1}(2)|"'] \arrow[luu, hook', "|\co_{\PP^1}(2)|", bend left] \arrow[ruu, hook, "|\co_{\PP^1}(2d-4)|"',  bend right] &             \\
      & \Gamma \arrow[u] \arrow[luuu,  rounded corners=10pt, 
    to path={(\tikztostart) -| (\tikztotarget) \tikztonodes}] \arrow[ruuu, rounded corners=10pt, 
    to path={(\tikztostart) -| (\tikztotarget) \tikztonodes}] & 
\end{tikzcd}
\]

\begin{proposition}
    \label{prop: P2 factoring RNC}
    Let $d \geq 3$ and $s = 2d-4$.
    Suppose that we are given $2d$ nondegenerate points $\Gamma$ in $\PP^s$ and a unique rational normal curve $C \to \PP^{s}$ through $\Gamma$, so that the embedding of $\Gamma$ in $\PP^s$ factors as $\Gamma \xrightarrow{\alpha}C \to \PP^s$. Then there exists a factorization $\Gamma \xrightarrow{\alpha} C \to \PP^2 \dashrightarrow \PP^s$, where the last map is induced by a linear series of $\co(d-2)$ of dimension $s+1 = 2d-3$, and there exists a $\frac{(d-2)(d-3)(2d-3)}{2}$-dimensional family of such factorizations.
\end{proposition}

\begin{proof}
    By Goppa duality for curves applied to the composition $\Gamma \xrightarrow{\alpha} C \to \PP^s$, we have a factorization $\Gamma \xrightarrow{\alpha} C \to \PP^2$ of the Gale transform, where the second map is a closed embedding describing a conic. Let $f$ be the equation of the conic in a parameterization of $\PP^2$. By Goppa duality for  $\PP^2$ (Corollary \ref{cor:complete intersection on P2}),  the original map factors in as $\Gamma \xrightarrow{\alpha} C \xrightarrow{\beta} \PP^2 \dashrightarrow \PP^s, $ where the last map is induced by a $(2d-3)$-dimensional linear series $W$ of $ \mathcal{O}(d-2)$ that is complementary to the subspace $H^0(\PP^2, \co(d-4)) \cdot f$.   Because the trivial intersection requirement is a proper closed condition, the moduli space for $W$ is a dense open subset of the family $G(2d-3, d(d-1)/2)$ of $(2d-3)$-dimensional subspaces of the vector space $H^0(\PP^2, \mathcal{O}(d-2))$.

    Now,  if two maps $\PP^2 \dashrightarrow \PP^s$ induced by the same linear series $W$ map $\Gamma$ to the same points in $\PP^s$, then they must be related by an automorphism of $\PP^s$ fixing the $2d$ points; this is possible if and only if the automorphism if the identity map.  Therefore, once we fix $W$, the map $\PP^2 \dashrightarrow \PP^s$ factoring the original map $C \to \PP^s$ is unique. Because $\dim G(2d-3, d(d-1)/2) = \frac{(d-2)(d-3)(2d-3)}{2}$, we obtain that the family of such factorizations of $\Gamma \to \PP^s$ is  $\frac{(d-2)(d-3)(2d-3)}{2}$-dimensional.
\end{proof}


An analogous statement holds for $(3,d)$-complete intersections, and the proof proceeds by the same method, by prescribing an elliptic curve of degree $3d-3$. We leave the details to the interested reader.


\begin{proposition}
    Let $d \geq 4$ and $s = 3d-4$.
    Suppose that we are given $3d$ nondegenerate points $\Gamma$ in $\PP^s$ and a unique elliptic curve $C \to \PP^{s}$ of degree $s+1 = 3(d-1)$ through $\Gamma$, so that the embedding of $\Gamma$ in $\PP^s$ factors as $\Gamma \xrightarrow{\alpha} C \to \PP^s$. Then there exists a factorization $\Gamma \xrightarrow{\alpha} C \to \PP^2 \dashrightarrow \PP^s$, where the last map is induced by an $(s+1)$-dimensional linear series of $\co(d-1)$, and there exists a $\frac{3(d-1)(d-2)(d-3)}{2}$-dimensional family of such maps $ \PP^2 \dashrightarrow \PP^s$ factoring $C \to \PP^s$.
\end{proposition}

\subsection{When $\Gamma$ is ``Almost'' a Complete Intersection}
\label{sec: almost complete intersection}
In previous sections, we focused on the case in which the configurations of points were realized as complete intersections on $\PP^2$. We note that in all but finitely many cases, the condition for the points to form a complete intersection on $\PP^2$ forces the points to lie on a special curve. In particular, the approach in Section \ref{sec: complete intersection} cannot be applied to points in general position in $\PP^2$ for a sufficiently large number of points in $\Gamma$ (i.e., $\gamma > 6$). For example, eight general points in $\PP^2$ do not form a complete intersection on $\PP^2$: a complete intersection of length 8 must be either a $(1,8)$-complete intersection or a $(2,4)$-complete intersection, but any curve passing through eight general points must be of degree at least three. Despite the eight general points failing to form a complete intersection on $\PP^2$, we note that the eight points are ``almost'' a complete intersection, in the sense that they lie on two degree-three curves: the eight points are one point away from being a $(3,3)$-complete intersection. To remove the excess point $p$ in the complete intersection, we use the method of blowup, which allows us to separate the two curves meeting at $p$ transversely in $\PP^2$. More precisely, by lifting the points and curves to the blowup surface $\Bl_p\PP^2$ under the blowup map $\Bl_p\PP^2 \to \PP^2$, we can achieve the eight general points as a complete intersection of two curves on $\Bl_p\PP^2$, both of class $3H-E$, where $H$ is the hyperplane class under the blowup map and $E$ is the exceptional divisor. Under this setting, we can safely apply Goppa duality for surfaces (Theorem \ref{thm: goppa duality for surfaces}) to obtain factorization of the Gale dual pair, through the blowup surface $\Bl_p\PP^2$. As a result, we obtain the following level diagram describing the duality for the case of eight general points in $\PP^4$:

    \[\begin{tikzcd}[column sep=1em, row sep=small]
\PP^2 = \PP H^0(\Bl_p\PP^2,\co(H))^\vee &  & \PP H^0(\Bl_p\PP^2,\co(2H-E))^\vee = \PP^4 \\
      & \Bl_p\PP^2 \arrow[lu, two heads, "|H|"] \arrow[ru, " |2H-E|"'] &             \\
      & C_i \arrow[u] &             \\
      & \Gamma \arrow[u] \arrow[luuu, 
       "\operatorname{res}(|H|) \subset |\co(H)|_\Gamma|", 
      rounded corners=10pt, 
    to path={(\tikztostart) -| (\tikztotarget) \tikztonodes}] \arrow[ruuu, 
     "\operatorname{res}(|2H-E|) \subset |\co(2H-E)|_\Gamma|"', 
    rounded corners=10pt, 
    to path={(\tikztostart) -| (\tikztotarget) \tikztonodes}] & 
\end{tikzcd}
\]
\bigskip

\noindent Generalizing this argument, in this section, we consider the case in which the points $\Gamma$ do not form a complete intersection on $\PP^2$, but become a complete intersection when we introduce additional points, which we denote by $R$. We then blowup $\PP^2$ at the set $R$ of the excess points, so that the lift of $\Gamma$ is realized as a complete intersection on the blown-up surface. 

\subsubsection{Case $|R| = 1$}
\label{sec: R = 1}

In this section, we make rigorous the discussion of eight general points $\Gamma$ in $\PP^2$ introduced in the previous section. Let $C_1$ and $C_2$ be two (degree-three) elliptic curves in $\PP^2$ passing through the eight general points $\Gamma$, and denote the excess intersection point by $p$. After blowing up $\PP^2$ at $p$, the points $\Gamma$ become the complete intersection of two elliptic curves $\widetilde{C}_1, \widetilde{C}_2$, both of class $3H-E$.

We intend to use Goppa duality for surfaces on $\Gamma \to \widetilde{C}_i \to \Bl_p\PP^2 \to \PP^2$, where the last map is the blowup map; the corresponding linear series is $|\co(H)|$ on $\Bl_p\PP^2$. Because the restriction map $H^0(\Bl_p\PP^2, \co(H)) \to H^0(\Gamma, \co(H)|_\Gamma)$ is injective, the requirement (ii) for Goppa duality for surfaces is satisfied.
We also need to check the vanishing of the first cohomology group $H^1$ of the following line bundles: $\mathcal{O}(H), \mathcal{O}(H) \otimes \mathcal{O}(-\widetilde{C}_1) , $ and  $ \mathcal{O}(H) \otimes \mathcal{O}(-\widetilde{C}_2)$. One can use the long exact sequence of cohomology for the closed subscheme exact sequence to compute $H^1$ of such line bundles.

\begin{lemma}
    \label{lem: vanishing of h1}
    The first cohomology groups $H^1$ of the line bundles $\mathcal{O}(H) $ and $\mathcal{O}(-2H+E) $  are trivial. 
\end{lemma}

\begin{proof}
     Because $\beta$ is the blowup map, we have the isomorphism $\beta_*\mathcal{O}_{\Bl_p\PP^2}  \cong \mathcal{O}_{\PP^2}$, so by the projection formula, we have an isomorphism  $\beta_*\beta^*\co(dH) \cong \co(dH)$ for any integer $d$.
    By the induced long exact sequence of cohomology for $$0 \to \mathcal{O}(-2H) \to  \mathcal{O}(-2H+E) \to  \mathcal{O}(-2H+E)|_E \to  0,$$ we have an exact sequence of cohomologies,
    \begin{align*}
     & H^1(\Bl_p\PP^2, \mathcal{O}(-2H)) \to  H^1(\Bl_p\PP^2, \mathcal{O}(-2H+E)) \to  H^1(E_p\PP^2, \mathcal{O}(-2H+E)|_E) \\ &\to H^2(\Bl_p\PP^2, \mathcal{O}(-2H))  ,
    \end{align*}
    which simplifies to  
    \begin{align*}
    0 \to  H^1(\Bl_p\PP^2, \mathcal{O}(-2H+E)) \to  H^0(E_p\PP^2, \mathcal{O}_E(-1)) \to 0.
    \end{align*}
    Because the exceptional divisor $E_p\PP^2$ is isomorphic to $\PP^1$, $H^0$ of $\co_E(-1)$ vanishes, hence $ H^1(\Bl_p\PP^2, \mathcal{O}(-2H+E)) = 0$.
\end{proof}
\noindent Therefore, Goppa duality for surfaces (Theorem \ref{thm: Goppa duality}) applies, and  the Gale transform $\Gamma \to \PP^4$ factors through the blowup surface, where the map  $\Bl_p\PP^2 \dashrightarrow \PP^4$ is induced by a 5-dimensional linear series of the Goppa dual line bundle $K_{\Bl_p\PP^2}(\widetilde{C}_1 + \widetilde{C}_2) \otimes \mathcal{O}(H)^\vee \cong \mathcal{O}(2H-E)$. By the projection formula, the global sections of $\co(2H-E)$ on $\Bl_p\PP^2$ correspond to quadratic sections on $\PP^2$ passing through $p$, which is five-dimensional. Therefore, the Goppa dual linear series must be the complete linear series $|2H-E|$. 


By starting with eight general points in $\PP^4$ and using Bézout's Theorem, we can reverse the discussion above and obtain a proof for the existence of a  factorization of the map of eight general points in $\PP^4$ through the surface $\Bl_p\PP^2$ as follows.
\begin{proposition}
    \label{prop: existence of factorization through BlpP2}
    Let $\Gamma$ be eight general points in $\PP^4$. Then there exists a factorization $\Gamma \to \Bl_p\PP^2 \xrightarrow{|2H-E|}  \PP^4$, where the map $\Gamma \to \Bl_p\PP^2$ is a closed embedding realizing $\Gamma$ as a complete intersection of two elliptic curves, both of class $3H-E$.
\end{proposition}

\begin{proof}
    The Gale transform of $\Gamma \to \PP^4$ is eight general points $\Gamma \hookrightarrow \PP^2$; by generality, there exists a unique pencil of cubic curves passing through the eight points in $\PP^2$. By Bézout's Theorem, there exists a unique point $p$ disjoint from $\Gamma$ such that the union $\Gamma \cup \{p\}$ is the intersection of any pair of distinct elements $C_1, C_2$ in the pencil; we may choose $C_1, C_2$ to be (smooth) elliptic curves. Consider the blowup of $\PP^2$ at $p$; the strict transform $\widetilde{C}_i$  of $C_i$  is an elliptic curve of class $3H-E$, and  $\Gamma$ is realized as the complete intersection of $\widetilde{C}_1$ and $\widetilde{C}_2$ on $\Bl_p\PP^2$. Hence we have a factorization of the Gale transform
    $\Gamma \to  \widetilde{C}_i \hookrightarrow \Bl_p\PP^2 \xrightarrow{\beta} \PP^2$, where $\beta$ is the blowup map.
    By Lemma \ref{lem: vanishing of h1},  Goppa duality for surfaces (Theorem \ref{thm: Goppa duality}) induces a factorization of the original map $\Gamma \to \PP^4$ as $\Gamma \to \widetilde{C}_i \hookrightarrow \Bl_p\PP^2 \dashrightarrow \PP^4,$
    where the last map is induced by the complete linear series of the Goppa dual line bundle
    $K_{\Bl_p\PP^2}(\widetilde{C}_1 + \widetilde{C}_2) \otimes \mathcal{O}(H)^\vee \cong \mathcal{O}(2H-E)$. 
    Therefore, eight general points in $\PP^4$ factor through a surface isomorphic to $\Bl_p\PP^2$ via the complete linear series $|2H-E|$.
    
\end{proof}

Since the surface $S$ is no longer the projective plane $\PP^2$, the uniqueness argument used in Lemma \ref{lem: (3,3)-CI uniqueness} where $S$ is $\PP^2$ does not directly apply. However, using the observation that the the automorphism group of $\Bl_p\PP^2$ can be regraded as the automorphism group of $\PP^2$ fixing a point, we obtain a similar unique factorization argument for Proposition \ref{prop: existence of factorization through BlpP2}.

\begin{proposition}
    \label{prop: unique fact through |2H-E|}
    Let $\Gamma$ be eight general points in $\PP^4$.  Then there exists a unique factorization $\Gamma \to \Bl_p\PP^2 \xrightarrow{|2H-E|}  \PP^4$, where $\Bl_p\PP^2$ achieves $\Gamma$ as a complete intersection of two curves of class $3H-E$, up to automorphism of $\Bl_p\PP^2$. 
\end{proposition}

\begin{proof}
    Given the closed embedding $\Gamma \to \PP^4$ of eight general points, by Proposition \ref{prop: existence of factorization through BlpP2}, we have a factorization $\Gamma \hookrightarrow \widetilde{C}_i \hookrightarrow \Bl_p\PP^2 \to \PP^4$, where the last map is induced by the complete linear series $|2H-E|$.
    Suppose that we have another factorization $\Gamma \xrightarrow{\alpha} \widetilde{C}_i' \xrightarrow{\pi'} \Bl_p\PP^2 \to \PP^4$, where $\widetilde{C}_1 \to \Bl_p\PP^2$ and $\widetilde{C}_2 \to \Bl_p\PP^2$ are the closed embeddings of curves of class $3H-E$ in $\Bl_p\PP^2$. By Goppa duality for surfaces (Theorem \ref{thm: Goppa duality}), its Gale transform $\Gamma \to \PP^4$ factors as $\Gamma \xrightarrow{\alpha'} \widetilde{C}_i' \xrightarrow{\pi'} \Bl_p\PP^2 \xrightarrow{\beta} \PP^2\xrightarrow{\kappa, \sim} \PP^2$. Let $p' = \kappa(p)$ be the image of the exceptional divisor of $\Bl_p\PP^2$ under $\Bl_p\PP^2 \xrightarrow{\kappa \circ \beta} \PP^2$. Then in the image in $\PP^2$, $\widetilde{C}_i'$ maps to a degree $3$ curve passing through $\Gamma$ and $p'$. By construction, $\Gamma$ in $\PP^2$ is the $(3,3)$-complete intersection of $\widetilde{C}_i$, so by uniqueness of the excess point of $(3,3)$-complete intersections passing through $\Gamma$, we must have $p = p'$, i.e., $\kappa$ must fix $p$. Hence, using the isomorphism $\operatorname{Aut}(\Bl_p\PP^2) \cong \operatorname{Aut}_p(\PP^2)$, $\kappa$ can be lifted canonically to an automorphism of $\Bl_p\PP^2$, so that $\Bl_p \PP^2 \xrightarrow{\beta} \PP^2 \xrightarrow{\kappa} \PP^2$ can be factored as $\Bl_p \PP^2  \xrightarrow{\widetilde{\kappa}, \cong} \Bl_p\PP^2 \xrightarrow{\beta} \PP^2$. Hence, by absorbing the automorphism of the blowup in $\beta'$, we may assume that the the Gale transform factors as $\Gamma \xrightarrow{\alpha'} \widetilde{C}_i' \xrightarrow{\pi'} \Bl_p\PP^2 \xrightarrow{\beta} \PP^2$. Because $\beta$ is injective away from the exceptional locus, the composition $\pi' \circ \alpha' : \Gamma \to \Bl_p\PP^2$ must coincide with the original factorization. Hence, we have two factorizations $\Gamma \hookrightarrow \Bl_p\PP^2 \to \PP^4$ of the same map, where the first map also coincides. Because the second map differs by an automorphism of $\PP^4$, and automorphisms of $\PP^4$ can fix at most six general points, by generality of the configuration of the points $\Gamma$ in $\PP^4$, the map $\Bl_p\PP^2 \to \PP^4$ must also coincide in the two factorizations. Therefore, the factorization is unique.
\end{proof}

\subsubsection{Case $|R| = 2$}
\label{sec: R = 2}
We continue to the case in which the points $\Gamma$ in $\PP^2$ are two points away from being a complete intersection. As a first example, we consider seven general points in $\PP^2$, which are two points away from being a $(3,3)$-complete intersection. We obtain a result and proof analogous to that of Proposition \ref{prop: existence of factorization through BlpP2}.
\begin{proposition}
    \label{prop: 7 general in P3}
    Given seven general points $\Gamma$ in $\PP^3$, there exists a map $\Bl_{\{p,q\}}\PP^2  \to \PP^3$ induced by the complete linear series $|2H-E|$ factoring the inclusion $\Gamma \hookrightarrow \PP^3$. 
\end{proposition}

\begin{proof}
    Starting with the map  $\Gamma \to \PP^3$ describing the seven general points, we obtain its Gale transform $\Gamma \to \PP^2$. Then there exist three linearly independent cubics in $\PP^2$ passing through $\Gamma$;
    choose two such curves $C_i$, which intersect at $\Gamma$ and at two additional points $p, q$, each with multiplicity one. By blowing up $\PP^2$ at $\{p,q\}$, we can regard $\Gamma$ as a complete intersection of two curves $\widetilde{C}_i$, each of class $3H-E$ on $\Bl_{\{p,q\}}\PP^2$. 
    By Serre duality, one can check that the required cohomology vanishing conditions are satisfied: $$h^1(\Bl_{\{p,q\}}\PP^2, \co(-2H + E)) = h^1(\Bl_{\{p,q\}}\PP^2,\co( -H)) = 0.$$ By Goppa duality for surfaces (Theorem \ref{thm: Goppa duality}), we obtain a factorization of the original map $\Gamma \to \PP^3$ through $\Bl_{\{p,q\}} \PP^2$ as
    $$ \Gamma \hookrightarrow \widetilde{C}_i \hookrightarrow \Bl_{\{p,q\}}\PP^2 \to \PP^3,$$
    where the last map is induced by the Goppa dual linear series $|2H-E|$. The key constructions are summarized in the following diagram.
    
    \[\begin{tikzcd}[column sep=0em, row sep=small]
\PP^2 = \PP H^0(\Bl_{\{p,q\}}\PP^2,\co(H))^\vee &  & \PP H^0(\Bl_{\{p,q\}}\PP^2,\co(2H-E))^\vee = \PP^3 \\
      & \Bl_{\{p,q\}}\PP^2 \arrow[lu, two heads, "|H|"] \arrow[ru, " |2H-E|"'] &             \\
      & C_i \arrow[u] &             \\
      & \Gamma \arrow[u] \arrow[luuu, 
       "\operatorname{res}(|H|) \subset |\co(H)|_\Gamma|", 
      rounded corners=10pt, 
    to path={(\tikztostart) -| (\tikztotarget) \tikztonodes}] \arrow[ruuu, 
     "\operatorname{res}(|2H-E|) \subset |\co(2H-E)|_\Gamma|"', 
    rounded corners=10pt, 
    to path={(\tikztostart) -| (\tikztotarget) \tikztonodes}] & 
\end{tikzcd}
\]

\end{proof}

\begin{remark}
    In contrast to the case $|R| = 1$, however, the uniqueness argument of the factorization breaks down when two points are blown up; this is because the geometry of the blowup introduces additional degrees of freedom in the choice of pencil of cubics. 
    Suppose we have another factorization $\Gamma \to  \Bl_{\{p,q\}}\PP^2 \xrightarrow{} \PP^3$ of the map $\Gamma \to \PP^3$, where $\Gamma$ is achieved as a complete intersection of two curves $D_i$, each of class $3H-E$. Then by Goppa duality  for surfaces, we obtain a factorization of the Gale transform $\Gamma \to D_i \to \Bl_{\{p,q\}}\PP^2 \xrightarrow{|H|} \PP^2$. We may write as $\Gamma \to D_i \to \Bl_{\{p,q\}}\PP^2 \xrightarrow{\beta} \PP^2 \xrightarrow{\kappa, \sim} \PP^2$ for some automorphism $\kappa$ of $\PP^2$. 
    By construction, $D_i$ maps to cubics through $\Gamma$ and the two points $p' = \kappa(p), q' = \kappa(q)$, which are the image of the exceptional divisor. However, there is no way to guarantee a relation between $\{p', q'\}$ and $\{p, q\}$; in general, they are induced by a different pencil of cubics through $\Gamma$.
\end{remark}


\subsubsection{Case $|R| = 3$} 
\label{sec: |R| = 3}
As we observed in Section \ref{sec: R = 2}, when there are more excess points in the complete intersection, we expect more complexity to be added to the factorization of the map of points in projective spaces. For example, in Section \ref{sec: R = 2}, the uniqueness argument used in Proposition \ref{prop: unique fact through |2H-E|} of Section \ref{sec: R = 1} does not hold for case $|R| = 2$ due to the emergence of an additional choice for curves describing the complete intersection, which affects the configuration of the excess points. 
Nevertheless, in some cases, we can control the complexity of the choices of points by requiring one of the curves making up the complete intersection to pass through the excess points in higher multiplicity; then fewer excess points are needed in describing the complete intersection.


To observe this approach in action, consider the case of nine general points in $\PP^5$. 
Coble proved that there exist exactly four Veronese surfaces passing through the nine general points $\Gamma$ \cite[Theorem~19]{Coble1922}; he first observed that the nine points determine a sextic elliptic curve $C$ as in Proposition \ref{prop: unique sectic elliptic cruve through 9 general points in P5}, then argued that any Veronese surface containing $\Gamma$ must contain $C$ using the uniqueness of $C$, and deduced that the embedding of $C$ in the Veronese surface corresponds to the square roots of a fixed degree six line bundle, four of which exist.  
In the following, we provide an alternative proof to this result by using Goppa duality for surfaces (Theorem \ref{thm: goppa duality for surfaces}), with a choice of three excess points to blow up.

The Gale transform of nine general points in $\PP^5$ is nine general points in $\PP^2$, so the minimal choice of degrees for the curves to induce a complete intersection containing $\Gamma$ is the pair $(3,4)$, whose complete intersection consists of 12 points. If we use the proof technique used in Proposition \ref{prop: 7 general in P3}, then $R$ is taken to be the three non-collinear excess points in the complete intersection. However, the family of quartic curves through nine points is five dimensional, so we do not expect a uniqueness argument on the factorization as we saw in Proposition \ref{prop: unique fact through |2H-E|} to hold.
Instead, consider a \emph{quintic} curve $C_2$ through $\Gamma$ in $\PP^2$; the complete intersection with the unique cubic curve $C_1$ through $\Gamma$ consists of 15 points, leading to six excess points. Now, suppose that that we require that $C_2$ pass through three points $R = \{p_1,p_2,p_3\}$ on $C_1$, each with multiplicity two. 
We can determine conditions on $R$ that guarantee the existence of a quintic through $\Gamma$ and having nodes at $R$ 
If such a triplet $R$ exists, then we can factor the Gale transform through the blowup surface $\Bl_R\PP^2$. The following Lemma shows that this is indeed possible.

\begin{lemma}
    \label{lem: nine general points through four veronese}
    Let $\Gamma$ be nine general points in $\PP^5$. Then there exist a factorization of the closed embedding $\Gamma \hookrightarrow \PP^5$ through the blowup of $\PP^2$ at three non-collinear points $R$, via $\Bl_R\PP^2 \to \PP^5$ induced by the complete linear series $|4H-2E|$.
\end{lemma}

\begin{proof}
    The map $\Gamma \to \PP^5$ of nine general points in $\PP^5$ corresponds to the map $\alpha : \Gamma \to \PP^2$ of nine general points in $\PP^2$ under Gale duality, which has a unique cubic elliptic curve $C_1$ passing through $\Gamma$ in $\PP^2$. We will pick three non-collinear points $R \subset C$ disjoint from $\Gamma$, so that there exists a quintic curve $C_2$ through $\Gamma$ with nodes at $R$ along $C_1$. For such a curve to exist, on $C_1$, the quintic section describing $C_2$ vanishes to second order at $R$, so the space $H^0(C_1, \co(5)|_{C_1} \otimes \co_{C_1}(-\Gamma - 2R))$ of such sections must be nontrivial. Since $C_1$ is an elliptic curve, this is possible if and only if the line bundle $\co_{C_1}(R)$ is a square root of the degree six line bundle $\mathcal{O}_{\PP^2}(5)|_{C_1}  \otimes \co_{C_1}(-\Gamma)$, so there are exactly four choices for the divisor $[R]$ satisfying this property. Conversely, each of these four choices gives rise to a quintic curve $C_2$ through $\Gamma$ with nodes at $R$ along $C_1$.
    
    Now, fix $R$ to be one of the four divisors and consider the blowup $\Bl_R\PP^2$ of $\PP^2$ at $R$. The curves $C_1$ and $C_2$ lift to curves $\widetilde{C}_1$ and $\widetilde{C}_2$ in $\Bl_R\PP^2$ of class $3H-E$ and $5H-2E$, respectively, and $\Gamma$ can be regarded as the complete intersection of $\widetilde{C}_1$ and $\widetilde{C}_2$ on $\Bl_R\PP^2$. One can compute using long exact sequences to deduce
    $$h^1(\Bl_R\PP^2, -2H+E) = h^1(\Bl_R\PP^2, -4H+2E) = 0,$$
    so Goppa duality for surfaces (Theorem \ref{thm: Goppa duality}) holds, and we obtain a factorization of the original map via a 6-dimensional linear series of $\co(4H-2E)$ on the blown-up surface $\Bl_R\PP^2$. Because the space $H^0(\Bl_R\PP^2, \co(4H-2E))$ of global sections correspond to quartic sections on $\PP^2$ with nodes at all points of $R$, its dimension is 6, so the map is induced by the complete linear series $|4H-2E|$. Therefore, we have a factorization
    $$\Gamma \xrightarrow{\widetilde{\alpha}} \Bl_R\PP^2\xrightarrow{|4H-2E|} \PP^5.$$    
    Hence we obtain the following level diagram of the duality.
    
    \[\begin{tikzcd}[column sep=0em, row sep=small]
\PP^2 = \PP H^0(\Bl_{R}\PP^2,\co(H))^\vee &  & \PP H^0(\Bl_{R}\PP^2,\co(4H-2E))^\vee = \PP^5 \\
      & \Bl_{R}\PP^2 \arrow[lu, two heads, "|H|"] \arrow[ru, " |4H-2E|"'] &             \\
      & C_i \arrow[u] &             \\
      & \Gamma \arrow[u] \arrow[luuu, 
       "\operatorname{res}(|H|) \subset |\co(H)|_\Gamma|", 
      rounded corners=10pt, 
    to path={(\tikztostart) -| (\tikztotarget) \tikztonodes}] \arrow[ruuu, 
     "\operatorname{res}(|4H-2E|) \subset |\co(4H-2E)|_\Gamma|"', 
    rounded corners=10pt, 
    to path={(\tikztostart) -| (\tikztotarget) \tikztonodes}] & 
\end{tikzcd}
\]

\end{proof}

We note that the map $\Bl_R\PP^2 \to \PP^5$ has as image a Veronese surface in $\PP^5$: this is due to the observation that we can decompose the map into a composition $\Bl_R\PP^2 \to \PP^2 \hookrightarrow \PP^5$, where the first map is induced by the complete linear series $|2H-E|$ and is surjective. Therefore, instead of factoring $\Gamma \to \PP^5$ through $\Bl_R\PP^2$, we may consider a further factorization through the Veronese surface $\PP^2$. It follows that the four choices of the divisor $[R]$ correspond to four different Veronese surfaces factoring $\Gamma \to \PP^5$, as summarized in the Proposition below.

\begin{proposition}[Coble {\cite[Theorem~19]{Coble1922}}]
    \label{prop: four veronese through 9 points}
    Let $\Gamma$ be nine general points in $\PP^5$. Then there exist precisely four factorizations of the closed embedding $\Gamma \hookrightarrow \PP^5$ through a Veronese surface, and the Veronese surface contains the unique sextic elliptic curve passing through $\Gamma$.
\end{proposition}
\begin{proof}
    From the proof of Lemma \ref{lem: nine general points through four veronese}, the nine general points $\Gamma \to \PP^5$ have a factorization through the blowup $\Bl_R\PP^2$ of $\PP^2$ at three points $R$ via the map $\Bl_R\PP^2 \to \PP^5$ induced by the complete linear series $|4H-2E|$, and there are four choices for the divisor $[R]$ on the unique elliptic curve $C_1$ passing through the nine general points of the Gale transform. For each such $R$, we obtain a quintic curve $C_2$ whose equation is given by the basis element of $H^0(C_1, \co(5)|_{C_1} \otimes \co_{C_1}(-\Gamma - 2R))$, satisfying the property that $C_2$ passes through $\Gamma$ and has nodes at each point of $R$.
    
    On the other hand, note that the line bundle $\co(2H-E)$ on $\Bl_R\PP^2$ is the unique square root of $\co(4H-2E)$, and because $h^0(\Bl_R\PP^2, 4H - 2E ) = 6$ and $h^0(\Bl_R\PP^2, 2H - E ) = 3$, we have a composition of two maps
    $$\Bl_R\PP^2\xrightarrow{|2H-E|} \PP^2 \xrightarrow{|\mathcal{O}(2)|} \PP^5$$
    whose associated line bundle is $4H-2E$. Because both maps are non-degenerate and well-defined, the linear series inducing the composition is the complete linear series $|4H-2E|$. Therefore, by Goppa duality for surfaces (Theorem \ref{thm: Goppa duality}), the map of nine general points $\Gamma \to \PP^5$ factors through the blowup surface $\Bl_R\PP^2$ as 
    $$\Gamma \to \Bl_R\PP^2 \xrightarrow{|2H-E|} \PP^2 \hookrightarrow \PP^5.$$
    Hence, the level diagram for the duality has the following further factorization.
    
    
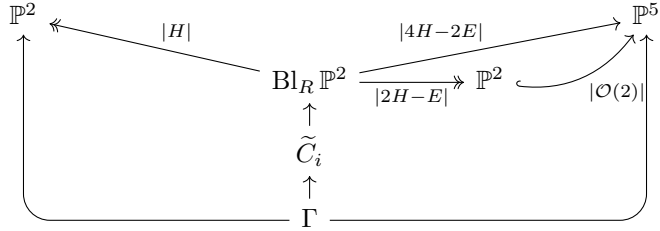
\[\begin{tikzcd}[column sep=4em, row sep=small]
\PP^2  & & & &  \PP^5 \\
      & & \Bl_{R}\PP^2 \arrow[llu, two heads, "|H|"'] \arrow[rru, " |4H-2E|"]  \arrow[r, "|2H-E|"', two heads] & \PP^2  \arrow[ru, hook, "|\co(2)|"', bend right]  &         \\
      & & \widetilde{C}_i \arrow[u] &  &           \\
      & &  \Gamma \arrow[u] \arrow[lluuu, 
      rounded corners=10pt, 
    to path={(\tikztostart) -| (\tikztotarget) \tikztonodes}] \arrow[rruuu, , 
    rounded corners=10pt, 
    to path={(\tikztostart) -| (\tikztotarget) \tikztonodes}] &  & 
\end{tikzcd}
\]

\bigskip

\noindent  We claim that the factorization through the Veronese surface is unique up to the choice of the divisor $[R]$ on the elliptic curve $C_1$. 
    Suppose that we have another factorization of the map $\Gamma \to \PP^5$,
    $$\Gamma \xrightarrow{\widetilde{\alpha}} \Bl_{R'}\PP^2\xrightarrow{|2H-E'|} \PP^2 \hookrightarrow \PP^5,$$
    where $R'$ describe the same divisor as $R$ on the elliptic curve $C_1$ in the blow-down $\PP^2$ of $\PP^2$, and $E'$ is the exceptional divisor of $\Bl_{R'}\PP^2$. By construction, there exists an elliptic curve $\widetilde{C}$ in $\Bl_R\PP^2$ of class $3H-E$ isomorphic to $C$, and the induced map $\widetilde{C} \to \PP^5$ is a closed embedding describing a sextic elliptic curve. By a similar argument imposed on the other factorization, we have another map $\widetilde{C}' \to \PP^5$ describing a sextic elliptic curve. 
    By Proposition \ref{prop: unique sectic elliptic cruve through 9 general points in P5}, the two maps must describe the same closed embedding of the elliptic curve. 
    It follows that the maps $\widetilde{C} \hookrightarrow \PP^2$ and  $\widetilde{C}' \hookrightarrow \PP^2$ into the Veronese surface are induced by the same line bundle of degree three because $R$ and $R'$ describe the same divisor on $C_1$. After a change of coordinates on $\PP^2$, we may assume that they describe the same map, so that we have two compositions $\widetilde{C} \hookrightarrow \PP^2 \hookrightarrow \PP^5$ of the same map, where the first maps also coincide. Because the composition describes the closed embedding of the unique sextic curve through $\Gamma$, it follows that the second map must also coincide, hence the factorization is unique if $R$ and $R'$ describe the same divisor on the curve $C$. 

    Conversely, if two factorizations come from different choice of the divisor $[R]$, then they describe distinct Veronese embeddings of $\PP^2$ in $\PP^5$ because the degree six maps from  the elliptic curve into $\PP^5$ are distinct. We conclude that there are exactly four Veronese surfaces passing through the nine general points in $\PP^5$.  
\end{proof}




    



\subsection{Limitations for $|R| > 3$}

In the examples above, we explicitly computed the dimension of the cohomology of line bundles on the blowup of $\PP^2$ to check the vanishing condition on the first cohomology groups $H^1$ required in using Goppa duality.
In general, as long as we can compute the dimensions of cohomology of line bundles involved in the computation, we can apply Goppa duality for surfaces to the cases with more points in the set $R$ of excess points. When the number of excess points is small, e.g. 1, 2, or 3, then we can use a linear change of coordinates to send the points in positions of our choice. However, when $R$ consists of more than three points, then the configuration of the points are relevant in the computation of cohomology dimensions, adding more complexity to applications of Goppa duality for surfaces.

\section{Relation to Prior Results and Potential Generalizations}
\label{sec: 4}
In this final section, we discuss the relation of Goppa duality for surfaces with prior results in the field, and also provide potential generalizations of Goppa duality for surfaces.

\subsection{Relation to prior results by  Deopurkar and Patel}

In \cite{DeopurkarPatel2024}, Anand Deopurkar and Anand Patel gave a proof for the count of Veronese surfaces in  passing through thirteen general points in $\PP^9$. Their proof idea can be interpreted using the technique of Goppa duality for surfaces described particularly in the proof of Proposition \ref{prop: four veronese through 9 points}. Consider thirteen general points $\Gamma$ in $\PP^9$; the Gale transform corresponds to thirteen general points in $\PP^2$.  On $\PP^2$, there exists a unique pencil of quartic curves passing through $\Gamma$, and there exists a triplet of points $R = \{p_1, p_2, p_3\}$ disjoint from $\Gamma$ such that any two distinct elements $C_1, C_2$ in the pencil intersect exactly at $\Gamma \cup R$. Take the simplistic guess to pick two quartic curves through $\Gamma$ to form an almost complete intersection. Then we can blow up $\PP^2$ at $R$ to lift $\Gamma$ to thirteen points on the blowup surface $\Bl_R\PP^2$, and lift $C_i$ to construct curves $\widetilde{C_i}$ of class $4H-E$, so that $\Gamma$ is a complete intersection on $\Bl_R\PP^2$. However, the Goppa dual line bundle is $\co(4H-E)$, and the dimension of the space of global sections is $h^0(\Bl_R\PP^2, \co(4H-E)) =  12,$
so the map to $\PP^9$ is rational in general and does not have an obvious factorization through $\PP^2$. 

Instead of choosing two quartic curves through $\Gamma$, consider choosing \emph{quintic} curves through $\Gamma$, which satisfy the additional condition that they must pass through three prescribed points $R = \{p_1, p_2, p_3\}$, each with multiplicity two. We wish conditions on $R$ which allow for a pencil of such quintics to exist. Using incidence correspondence, one can check that there are only finitely many configurations of $R$ which satisfy this condition.   Now, suppose we fix a configuration $R$ that satisfy this property; then we can blow up $\PP^2$ at $R$ to obtain two curves $\widetilde{C}_i$ of class $5H-2E$ that pass through the lift of $\Gamma$ in $\Bl_R\PP^2$. The Goppa dual line bundle is then $\co(6H-3E)$, with $h^0(\Bl_R\PP^2, \co(6H-3E)) = 10$, so the induced map from $\Bl_R\PP^2$ to $\PP^9$ is induced by the complete linear series $|6H-3E|$. Moreover, the map factors through $\PP^2$ via $\co(3)$, so the thirteen points $\Gamma$ lie on a Veronese surface $\PP^2$ in $\PP^9$. The key constructions are summarized in the following diagram.

    
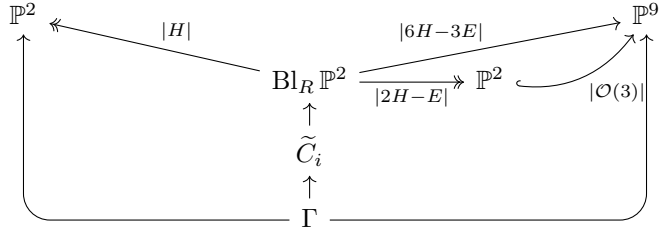
\[\begin{tikzcd}[column sep=4em, row sep=small]
\PP^2  & & & &  \PP^9 \\
      & & \Bl_{R}\PP^2 \arrow[llu, two heads, "|H|"'] \arrow[rru, " |6H-3E|"]  \arrow[r, "|2H-E|"', two heads] & \PP^2  \arrow[ru, hook, "|\co(3)|"', bend right]  &         \\
      & & \widetilde{C}_i \arrow[u] &  &           \\
      & &  \Gamma \arrow[u] \arrow[lluuu, 
      rounded corners=10pt, 
    to path={(\tikztostart) -| (\tikztotarget) \tikztonodes}] \arrow[rruuu, , 
    rounded corners=10pt, 
    to path={(\tikztostart) -| (\tikztotarget) \tikztonodes}] &  & 
\end{tikzcd}
\]

Conversely, if we have thirteen general points $\Gamma$ lying on a Veronese surface in $\PP^9$, then on the surface $\PP^2$, we have a pencil of quartic curves $D_i$ through $\Gamma$, whose elements pass through three excess points $T = \{q_1,q_2, q_3\}$. We can blow up $\PP^2$ at these points to achieve $\Gamma$ as a complete intersection on two curves $\widetilde{D}_i$ of class $4H-E$ on  $\Bl_R\PP^2$. By Goppa duality for surfaces (Theorem \ref{thm: goppa duality for surfaces}), with the choice $L = \co(3H)$ and $W = |L|$, we obtain that the Goppa dual line bundle is $\co(2H-E)$ and the Goppa dual linear series is the complete linear series $|2H-E|$. Under the map $\Bl_T\PP^2 \xrightarrow{|2H-E|} \PP^2$, the curves $\widetilde{D}_i$ map to quintic curves with nodes precisely at $R' = \{p_1', p_2', p_3'\}$.  Hence we have another level diagram given below.

    
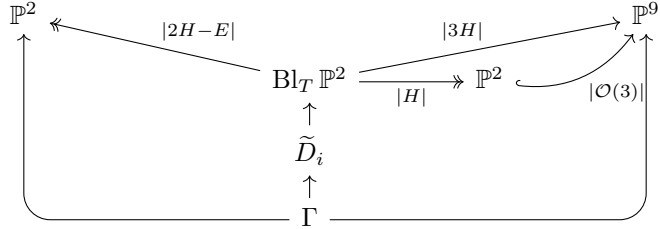
\[\begin{tikzcd}[column sep=4em, row sep=small]
\PP^2  & & & &  \PP^9 \\
      & & \Bl_{T}\PP^2 \arrow[llu, two heads, "|2H-E|"'] \arrow[rru, " |3H|"]  \arrow[r, "|H|"', two heads] & \PP^2  \arrow[ru, hook, "|\co(3)|"', bend right]  &         \\
      & & \widetilde{D}_i \arrow[u] &  &           \\
      & &  \Gamma \arrow[u] \arrow[lluuu, 
      rounded corners=10pt, 
    to path={(\tikztostart) -| (\tikztotarget) \tikztonodes}] \arrow[rruuu, , 
    rounded corners=10pt, 
    to path={(\tikztostart) -| (\tikztotarget) \tikztonodes}] &  & 
\end{tikzcd}
\]

\noindent Using this observation, counting the number of Veronese surfaces through thirteen general points in $\PP^9$ roughly reduces to counting all possible choices for the triplet $R$, and computing the number of factorizations of the Gale transform associated with each choice of $R$. Deopurkar and Patel perform the rigorous computation of this number in their paper \cite{DeopurkarPatel2024}, by parametrizing the triplet of points in terms of \emph{complete triangles}, which is a $\operatorname{PGL}(3)$-equivariant modification of the Hilbert scheme $\operatorname{Hilb}_3\PP^2$ of length-three subschemes in $\PP^2$.

\subsection{Potential generalizations}

In this paper, we generalized the theory of classical Goppa duality to the case where the points were achieved as a complete intersection on a surface with a map to projective spaces. We note that Theorem \ref{thm: 3.1} by Eisenbud-Popescu on the existence of an exact sequence involving the spaces of global sections and the dual spaces holds for any locally free resolution of a line bundle on $\Gamma$. In particular, we can derive Goppa duality on higher-dimensional schemes, which are induced by projective vector bundles $\PP (\mathcal{E}_0)$ over $B$. For example, if $B$ is a curve and $\mathcal{E}_0$ is of rank 2, then we have a ruled surface $\PP\mathcal{E}_0$ over $B$, and $\Gamma$ can be regarded as points on $\PP\mathcal{E}_0$. Moreover, this method allows the map $\Gamma \to B$ to be not necessarily a closed embedding but a finite morphism. This method is particularly useful when we fail to find a simpler curve passing through the points, but some ruled surface over some curve contains the points.

We can also generalize Goppa duality for surfaces to higher-dimensional schemes. For example, Theorem \ref{thm: goppa duality for surfaces} can be generalized to arbitrary zero-dimensional complete intersections of a Gorenstein scheme $B$. As the dimension of the base scheme $B$ increases, we have more vanishing constraint of cohomology as stated in Theorem  \ref{thm: 3.1}; hence this generalization would be of good use particularly when the base scheme has simple cohomological properties, such as projective spaces $\PP^c$, which guarantee the vanishing of intermediate cohomologies on line bundles.

\bibliographystyle{amsplain}
\bibliography{refs} 
\end{document}